\documentclass{amsart}

\usepackage{amsfonts, amssymb, amsmath, amsthm}                               
\usepackage{geometry}                                                         
\usepackage{setspace}                                                         
\usepackage{url}                                                              
\usepackage{color}
\usepackage{slashed}
\usepackage{comment}
\usepackage{hyperref}

\newtheorem{theorem}{Theorem}

\newtheorem{lemma}[theorem]{Lemma}
\newtheorem{proposition}[theorem]{Proposition}

\newtheorem{assumption}[theorem]{Assumption}
\newtheorem{remark}[theorem]{Remark}

\numberwithin{equation}{section}                                              
\numberwithin{theorem}{section}                                               
\numberwithin{figure}{section}                                                

\allowdisplaybreaks
\begin{document}
\title[Hyperbolic systems with non-diagonalisable principal part and variable multiplicities]{Hyperbolic systems with non-diagonalisable principal part and variable multiplicities, IV: time-dependent characteristics, representation of solutions and regularity}

\author{Michael Ruzhansky}
\address{Michael Ruzhansky: Department of Mathematics: Analysis, Logic and Discrete Mathematics, Ghent University, Krijgslaan 281, Building S8, B 9000 Ghent, Belgium; and School of Mathematical Sciences, Queen Mary University of London, Mile End Road, London E1 4NS, United Kingdom}
\email{Michael.Ruzhansky@UGent.be; m.ruzhansky@qmul.ac.uk}

\author{Bolys Sabitbek}
\address{Bolys Sabitbek:  Department of Mathematics: Analysis, Logic and Discrete Mathematics, Ghent University, Krijgslaan 281, Building S8, B 9000 Ghent, Belgium; and Institute of Mathematics and Mathematical Modeling, Almaty, Kazakhstan}
\email{sabitbek.bolys@ugent.be}

\begin{abstract}
We consider first-order linear hyperbolic systems with time-dependent pseudo-differential principal symbols and characteristic roots of variable multiplicity. Building on the time-independent work of Kamotski and the first author, we introduce a finite-type condition along the lifted bicharacteristic flow that accounts for the contact of time-dependent characteristic roots. Under this condition, we construct a representation of solutions using time-dependent Fourier integral operators and a Picard series whose terms become increasingly smooth. This leads to an asymptotic expansion in smoothness and local \(L^p\) regularity estimates with the expected loss of derivatives. We also apply these results to higher-order scalar weakly hyperbolic equations with involutive characteristics and time-dependent principal symbols. Finally, we consider upper triangular systems with possibly non-diagonalisable principal parts and derive a representation formula involving integrated Fourier integral operators. This extends the earlier results of Garetto, Jäh and the first author to time-dependent principal parts.
\end{abstract}
\maketitle

\section{Introduction}
The present paper has two objectives. First, we extend the representation formulas and regularity results of Kamotski and Ruzhansky \cite{KR2007} to hyperbolic
systems with time-dependent principal symbols and characteristic roots of variable multiplicity. Second, we continue the study of hyperbolic systems with non-diagonalisable principal parts and variable multiplicities initiated in \cite{GJR18} and developed further in \cite{GJR20,GS24}. In this fourth part of
the series, we establish representation formulas and regularity estimates for such systems with time-dependent principal symbols.
\subsection{Problem setting}
We consider the Cauchy problem for first-order hyperbolic systems
\begin{equation}\label{eq:CP}
\begin{cases}
D_t u = P(t,x,D_x)\,u, & (t,x)\in[0,T]\times\mathbb R^n,\\
u(0,x) = u_0(x),
\end{cases}
\end{equation}
where $D_t=-i\partial_t$ and $P$ is an $m\times m$ matrix of classical pseudo-differential operators of order one, depending smoothly on $t$. The principal symbol of $P$ has smooth real eigenvalues $a_1,\dots,a_m$, which depend on $(t,x,\xi)$ and may coincide with variable multiplicities. We study two classes of such systems:
\begin{enumerate}
  \item systems whose principal symbol is microlocally smoothly diagonalisable (Sections~\ref{sec:picard}--\ref{sec:principle});
  \item upper triangular systems, which need not be diagonalisable (Section~\ref{sec:triangular-representation}).
\end{enumerate}
For the first class, we extend the approach of \cite{KR2007} by introducing a finite-type condition adapted to the time-dependent characteristic roots. For the second, we extend the representation theory of \cite{GJR20} to time-dependent principal symbols.

\subsection{Background and motivation} 

\subsubsection{Strictly hyperbolic systems and constant multiplicity}

For strictly hyperbolic systems, the principal symbol can be diagonalised, and the solution of \eqref{eq:CP} can be represented as a sum of Fourier integral operators acting on the initial data \cite{H71,Du96}. The calculus of Fourier integral operators then provides a well-established framework for studying regularity and the propagation of singularities \cite{DH72}, while the spectral asymptotics of $P$ are described by the Weyl formula \cite{H68,DG75}.
When the characteristic roots have constant multiplicity, the principal symbol can be reduced to block-diagonal form, allowing much of the strictly hyperbolic theory to be carried over to this setting. We refer to the works of Taylor \cite{T81}, Kumano-go \cite{K81}, and Kumano-go, Taniguchi and Tozaki \cite{KTT78}. Further accounts can be found in Ivrii \cite{Iv93} and, for analytic coefficients, in Nishitani \cite{Ni14}.

\subsubsection{Variable multiplicity and finite-type condition} 
When eigenvalues of a system have variable multiplicity, the block-diagonal structure breaks down, making direct diagonalisation impossible. Kucherenko \cite{Kur74} addressed this difficulty by constructing formal asymptotic solutions in this setting. 
The first step toward a rigorous theory came from Rozenblum \cite{Roz80}, who built the fundamental solution using the Picard series, a method of successive approximation. He proved that this solution admits a well-controlled asymptotic expansion with respect to smoothness. His approach relied on a transversality condition: wherever two eigenvalues $a_j$ and $a_k$ coincide at a point $(x,\xi)$, their Poisson bracket must satisfy 
\begin{equation*}
    \{ a_j, a_k \}(x,\xi) \neq 0.
\end{equation*}
Kamotski and the first author \cite{KR2007} later extended this theory by relaxing that requirement. Instead of full transversality, they imposed the weaker finite-type condition: for some integer $\nu \geq 1$, the iterated Poisson brackets satisfy 
\begin{equation*}
   \{a_j, \{ a_j, \ldots, \{ a_j, a_k \} \} \ldots \}^{(\nu)}(x,\xi) \neq 0. 
\end{equation*}
This condition permits so-called involutive characteristics, where the Poisson bracket itself vanishes. Under these more general geometric conditions on the characteristics, they showed that solutions are represented as  a Picard series whose terms are iterated Fourier integral operators. This series serves as an asymptotic expansion with respect to smoothness, and it yields a precise description of propagation of singularities together with sharp regularity properties of solutions. These results were applied in two directions: to establish regularity estimates for scalar weakly hyperbolic equations with involutive characteristics, and to derive the first and second terms of the spectral asymptotics for the corresponding self-adjoint elliptic systems.  

Building on this framework, Garetto, J\"ah, and the first author proved well-posedness of the Cauchy problem in anisotropic Sobolev spaces for non-diagonalisable hyperbolic systems with variable multiplicities \cite{GJR18}, and established a representation formula for the solution in terms of Fourier integral operators \cite{GJR20}, which in turn yields results on regularity and propagation of singularities; see \cite{GS24} for the case of singular coefficients and \cite{GS25} for the higher order hyperbolic equations with multiplicities. These results extend those of \cite{KR2007} to the non-diagonalisable case, although the representation formula of \cite{GJR20} still assumes that the principal symbol is independent of time.

The time-independent principal symbol is a restrictive assumption. For example, in many physical systems such as in elasticity, in electromagnetism in inhomogeneous media, and in geophysical wave propagation, the coefficients depend on time. Without this assumption, we can extend the theory to a wider class of problems. 
Therefore, this paper focuses on the extension of the framework of \cite{KR2007} to the principal symbol depending on time and spatial variables.  This is not a direct generalisation: when the time-dependent principal symbol changes the geometry of the characteristics and the structure of the Hamiltonian flow, and the Picard series must be constructed differently. With these modifications, we obtain the representation formulas, regularity estimates and propagation results in the general setting. Moreover, the same approach works also for non-diagonalisable upper triangular systems.

\subsection{The geometric condition}\label{subsec:condC}
In \cite{KR2007}, using iterated Poisson brackets, they measure the order of contact between characteristic surfaces as follows
\begin{equation*}
    H_{ a_j }^{ \nu } a_k = \underbrace{ \{ a_j, \{ a_j, \ldots, \{ a_j, a_k \} \ldots \} \} }_{ \nu \text{ times}}.
\end{equation*}
When $a_j$ varies with time, the bicharacteristic flow gains an explicit time component, so the spatial Poisson bracket alone fails to represent the contact geometry. We suggest the operator as a natural substitute
\begin{align*}
\mathcal{D}_j f := \partial_t f + H_{a_j}(f), \quad j = 1,\ldots,m,
\end{align*}
where
\begin{align*}
H_g(f) = \{g, f\} = \sum_{k=1}^{n} \left( \frac{\partial g}{\partial \xi_k}
\frac{\partial f}{\partial x_k} - \frac{\partial g}{\partial x_k}
\frac{\partial f}{\partial \xi_k} \right)
\end{align*}
denotes the Poisson bracket in the spatial variables. The integral curves of $\mathcal{D}_j$ are the lifted bicharacteristics of $a_j$ in the extended phase space $[0,T] \times T^*\mathbb{R}^n$:
\begin{align*}
\dot{t} = 1, \quad \dot{x} = \partial_{\xi} a_j(t,x,\xi), \quad \dot{\xi} = -\partial_x a_j(t,x,\xi).
\end{align*} 
Differentiation along these curves, rather than along the purely spatial Hamiltonian flow, is what the time-dependent setting demands.

We assume that, at points where two non-identical characteristic roots coincide, $a_j = a_k$, the lifted bicharacteristics generated by $a_j$ meet the coincidence set $\{a_j = a_k\}$ with finite order of contact. Equivalently, there exists a uniform integer $M\geq 1$ such that, at every such point, $\mathcal{D}_j^{ \nu }( a_j - a_k) \neq 0$ for some $1 \leq \nu \leq M$. This leads to the following finite-type contact condition along the lifted bicharacteristic flow.

\medskip
\noindent\textbf{Condition~C.}\label{cond:C}\;
\textit{The eigenvalues $a_1, \ldots ,a_m $ of the principal symbol of $P$ are smooth and real, and there exists an integer $M\geq 1$ such that for every $(t,x,\xi) \in [0,T] \times T^* \mathbb{R}^n$ with $\xi \neq 0$ and every pair $j \neq k$, if $a_j(t,x,\xi) = a_k(t, x, \xi)$ and $a_j \not \equiv a_k$ near $(t, x, \xi)$, then there exists an integer $1 \leq \nu \leq M$ for which}
\begin{equation}\label{eq:ConditionC}
\mathcal{D}_j^{\nu}(a_j - a_k)(t,x,\xi)
\;=\; \bigl( \partial_t + H_{a_j} \bigr)^{\nu - 1}\bigl( \partial_t (a_j - a_k) + H_{a_j}(a_j - a_k)\bigr)(t,x,\xi)\;\neq\;0.
\end{equation}

\medskip
\noindent
For the time-independent eigenvalues, the operator $\mathcal{D}_j$ reduces to $H_{a_j}$, and Condition~C recovers the autonomous condition of \cite{KR2007}:
\begin{equation*}
\mathcal{D}_j^{ \nu }(a_j - a_k)(x, \xi)
\; = \; H_{ a_j }^{ \nu }( a_j - a_k)(x, \xi)
\;=\; - \overbrace{
\{ a_j, \{ a_j, \ldots, \{ a_j, a_k\} \} \ldots \}
}^{ \nu \text{ times}}(x, \xi)\;\neq\;0,
\quad 1\leq \nu \leq M.
\end{equation*}
In this context, Condition~C is the natural time-dependent generalisation of the finite-type condition introduced in \cite{KR2007}, serving as the geometric foundation for the subsequent analysis.

\subsection{Main results}\label{subsec:main-results}
Before stating the precise theorems, we describe our results informally. For diagonalisable systems, the solution of \eqref{eq:CP} can be written as $u=EVu_0$, where $E$ is a propagator built from time-dependent Fourier integral operators associated with the lifted bicharacteristic flows of the eigenvalues, and $V$ is a Picard series. Under Condition~C, the terms of this series become progressively smoother, so the Picard series is an asymptotic expansion in smoothness. This smoothing property yields $L^p$ regularity for solutions of \eqref{eq:CP}, as well as $L^p$ estimates for higher-order scalar weakly hyperbolic equations with involutive characteristics, for all $1<p<\infty$. Finally, for upper triangular systems, which need not be diagonalisable, we obtain a representation formula for the solution in terms of integrated Fourier integral operators.

For diagonalisable systems, we decompose $P = A + B$, where $A$ is a smooth time-dependent family of pseudo-differential operators of order one whose principal symbol is microlocally smoothly diagonalisable, and $B$ collects the lower-order terms as a smooth time-dependent family of pseudo-differential operators of order zero. The solution of~\eqref{eq:CP} is then written as $u(t)= E(t)V(t)u_0$, where $E(t)$ is the propagator of the diagonal principal part and $V(t) = I+\sum_{\ell\ge 1}Q_\ell(t)$ is the Picard series correcting for the off-diagonal contributions.  Both $E$ and $V$ are constructed precisely in Section~\ref{sec:prelim}, with the Picard series given in~\eqref{eq:Picard}--\eqref{eq:Ql}.

Our first main result shows that the iterated Picard terms $Q_{\ell}$ gain smoothness as $\ell$ increases. Consequently, the series is an asymptotic expansion in smoothness, and the regularity results that follow are derived from this property.

\begin{theorem}[Smoothing of the Picard series]\label{thm:smoothing}
Assume that $A(t,x,\xi)$ is microlocally smoothly diagonalisable and that its eigenvalues $a_{j}(t,x,\xi)$ satisfy Condition~C. Then there exists a sequence $0 \leq N(\ell)\to\infty$ such that each operator $Q_{\ell}(t)$ defined in~\eqref{eq:Ql} extends to a bounded map
\[
Q_{\ell}(t)\colon L^{2}_{\mathrm{comp}}(\mathbb{R}^n)\;\longrightarrow\; H^{N(\ell)}_{\mathrm{loc}}(\mathbb{R}^n),\qquad t\in[0,T].
\]
In particular, the Picard series~\eqref{eq:Picard} is an asymptotic expansion with respect to smoothness.
\end{theorem}

The smoothing property of Theorem~\ref{thm:smoothing}, combined with the abstract regularity principle of Theorem~\ref{thm:3.1}, transfers regularity from the data to the solution. This yields the following $L^p$ estimates for solutions of~\eqref{eq:CP}.

\begin{theorem}[$L^{p}$-regularity for systems]\label{thm:Lp-reg}
Let $1<p<\infty$ and set $\alpha:=(n-1)|1/p-1/2|$. Assume that $A$ is microlocally smoothly diagonalisable with smooth eigenspaces and real eigenvalues $ a_j(t, x, \xi)$ satisfying Condition~C. Then, for every compactly supported initial datum $u_0 \in L^p_{\alpha} \cap L^2_{\mathrm{comp}}$, the solution $u$ of~\eqref{eq:CP} satisfies $u(t,\cdot)\in L^p_{\mathrm{loc}}$ for all $t\in[0,T]$, and there is a constant $C_T>0$ such that
\begin{equation}\label{eq:2.4}
\sup_{0 \leq t \leq T}\|u(t,\cdot)\|_{L^p_{\mathrm{loc}}} \;\leq \; C_T\,\|u_0\|_{L^p_{\alpha}}.
\end{equation}
\end{theorem}

Theorem~\ref{thm:Lp-reg} is stated for first-order systems, but its scope is broader: by reducing a higher-order scalar equation to a first-order system, the same regularity theory applies to scalar weakly hyperbolic equations with involutive characteristics. We make this explicit for equations of the form 
\begin{equation}\label{eq:scalar-L}
L \;=\; \Lambda_{1}\cdots\Lambda_{m}\;+\;\sum_{k<m}\sum_{j_{1},\ldots,j_{k}} b_{j_{1} \ldots j_{k}}(t,x,D_{x})\, \Lambda_{j_{1}}\cdots\Lambda_{j_{k}} \;+\; c(t,x,D_{x}),
\end{equation}
where each factor $\Lambda_{j} = D_{t} + \lambda_{j}(t,x,D_{x})$ has $\lambda_{j}$ of order one, and the coefficients $b_{j_{1}\ldots j_{k}}$ and $ c$ are smooth time-dependent families of pseudo-differential operators of order zero. We denote the principal symbol of $\lambda_{j}$ by $a_{j}(t,x,\xi)$. The operator $L$ is said to have  involutive characteristics if
\begin{equation*}
[\Lambda_{j},\Lambda_{k}]=\alpha_{jk}\Lambda_{j}+\beta_{jk}\Lambda_{k}+\gamma_{jk}
\end{equation*}
where $\alpha_{jk},\beta_{jk}$, and $\gamma_{jk}$ are smooth time-dependent families of pseudo-differential operators of order zero.

\begin{theorem}[$L^{p}$-regularity for scalar weakly hyperbolic equations]\label{thm:weak-eq}
Let $1 < p < \infty$, $ s\in \mathbb{R}$, and $\alpha = (n-1)|1/p-1/2|$. Suppose the principal symbols $a_j$ of $\lambda_j$ satisfy Condition~C. Let $u$ be the solution to the Cauchy problem
\begin{equation}\label{eq:2.6}
Lu = 0, \qquad \partial_t^j u(0,x) = g_j(x), \quad 0 \leq j \leq m-1,
\end{equation}
with compactly supported Cauchy data $g_j \in L^p_{\alpha + s + m - 1 - j} \cap L^2_{\mathrm{comp}}$. Then $u(t,\cdot)\in(L^p_s)_{\mathrm{loc}}$ for all $ t \in [0,T]$, and
\begin{equation}
\sup_{ t \in [0,T]} \bigl\| \partial^J u(t,\cdot) \bigr\|_{(L^p_s)_{\mathrm{loc}}} \;\leq\; C \sum_{j=0}^{m-1}\|g_j\|_{L^p_{\alpha + s + m-1-j}},
\end{equation}
where $\partial^J=\Lambda_{j_1}\cdots\Lambda_{j_k}$ ranges over all ordered tuples $(j_1, \ldots ,j_k) \in \{1, \ldots, m\}^k$ with $k \leq m-1$.
\end{theorem}

The results above concern diagonalisable systems. Our final result shows that the time-dependent framework also applies to non-diagonalisable systems. Following \cite{GJR18, GJR20}, we consider
\begin{equation}\label{eq:CP-triangular}
D_tu=A(t,x,D_x)u+B(t,x,D_x)u+f,\qquad u(0)=u_0,
\end{equation}
where $A$ is an upper triangular matrix of pseudo-differential operators of order one and $B$ is a matrix of pseudo-differential operators of order zero.

\begin{theorem}[Informal version of Theorem~\ref{rep:thm-main}]\label{intro:thm-main}
Assume that the lower order terms $B$ satisfy a Levi-type condition and that the diagonal symbols $a_1,\dots,a_m$ of $A$ satisfy Condition~C.
Let $u$ be the solution of \eqref{eq:CP-triangular} in anisotropic Sobolev spaces given by \cite{GJR18}, and denote by $u_l^0$ and $f_l$ the components of $u_0$ and $f$. Then, for every $N\in\mathbb N$,
\begin{equation*}
    u_j(t) = \sum_{l=1}^{m} \Bigl[ \bigl( \mathcal H_{j,l}(t) + R_{j,l}(t)\bigr) u_l^0
+ \bigl( \mathcal K_{j,l}(t) + S_{j,l}(t)\bigr)f_l\Bigr], \qquad j=1,\dots,m,
\end{equation*}
where $\mathcal H_{j,l}$ and $\mathcal K_{j,l}$ are integrated Fourier integral operators of order $l-j$, and the remainders $R_{j,l}$, $S_{j,l}$ gain $N$ derivatives with respect to them.
\end{theorem}

This extends the representation formula of \cite{GJR20} to time-dependent principal parts, and completes the study initiated in \cite{GJR18, GJR20} and continued in \cite{GS24}.

\subsection{Strategy of proof}\label{subsec:strategy}
Following \cite{KR2007}, we do not apply the multiphase calculus of Kumano-go \cite{K81} and Kumano-go--Taniguchi--Tozaki \cite{KTT78} to $Q_\ell$.
Recall that the Picard term $Q_\ell$ is an iterated time-ordered integral over the simplex $\Delta_\ell$, summed over all admissible mode chains $(j_1, \ldots, j_{\ell + 1})$ with $j_k \neq j_{k+1}$, whose integrand alternates between the Fourier integral propagators $E_{j_k}(t_{k-1}, t_k)$ and zero-order coupling terms. In the strictly hyperbolic setting, multiphase calculus collapses this iterated product into a single oscillatory integral with a phase assembled from the separate mode phases. Under variable multiplicity, however, eigenvalues coalesce, the individual phases fail to stay distinct, and the composite phase degenerates. Instead, we handle the iterated form of $Q_\ell$ directly: Egorov's theorem decouples the amplitude from the phase, a Hamilton--Jacobi system determines the composite phase, and the coincidence set is controlled via the finite-order contact hypothesis in Condition~C.

Under the substitution $u=EV u_0$, where $E$ is the propagator of the diagonal principal part, the Cauchy problem \eqref{eq:CP} reduces to an equation for $V$ whose $\ell$-th Picard term takes the form
\begin{equation*}
Q_\ell \;=\; \int_0^t \! \int_0^{t_1} \! \cdots \! \int_0^{t_{\ell-1}} \mathcal{B}(\bar t)\, H(\bar t)\,dt_\ell \cdots dt_1,
\end{equation*}
where $H(\bar t)$ is a composition of $\ell+1$ propagators $E_j$  and $\mathcal{B}(\bar t)$ is a smooth time-dependent family of pseudo-differential operators of order zero. The smoothing estimate for $Q_{\ell}$ is obtained in three stages: an amplitude-phase separation (Stage I), the construction of a phase function via a non-autonomous Hamilton--Jacobi system (Stage II), and a three-region decomposition of the resulting oscillatory integral (Stage III).

\subsection{Plan of the paper}\label{subsec:plan}
The remainder of the paper is organised as follows. Section~\ref{sec:prelim} fixes notation and recalls the microlocal diagonalisation that underlies the Picard construction. Section~\ref{sec:picard} develops the Picard analysis and proves Theorem~\ref{thm:smoothing}. Section~\ref{sec:principle} formulates an abstract regularity principle for first-order Cauchy problems and combines it with Theorem~\ref{thm:smoothing} to establish Theorems~\ref{thm:Lp-reg} and~\ref{thm:weak-eq}. Section~\ref{sec:triangular-representation} proves the representation formula of Theorem~\ref{rep:thm-main} for time-dependent upper triangular systems. The Appendix collects five technical results used throughout the paper: a lifting lemma, Egorov's theorem, an $L^2\to H^s$ Fourier integral operator estimate, the measure estimate under Condition~C, and an interpolation lemma.

\section{Preliminaries}\label{sec:prelim}

\subsection{Notation and function spaces}\label{subsec:notation}
We write $D_x = -i\partial_x$ and $\langle \xi \rangle = (1 + |\xi|^2)^{1/2}$ for $\xi \in \mathbb{R}^n$. For $s \in \mathbb{R}$, the Sobolev space $H^s(\mathbb{R}^n)$ consists of all tempered distributions $u \in \mathcal{S}' (\mathbb{R}^n)$ such that
\begin{equation*}
\|u\|_{H^s} = \|\langle \xi \rangle^s \widehat{u}(\xi)\|_{L^2} < \infty,
\end{equation*}
where $\widehat{u}$ denotes the Fourier transform of $u$. For $s \in \mathbb{R}$  and $1 < p < \infty$, the Sobolev space $L^p_s(\mathbb{R}^n)$ consists of all tempered distributions $u \in \mathcal{S}'(\mathbb{R}^n)$ such that $(I - \Delta)^{s/2} u \in L^p(\mathbb{R}^n)$, with norm
\begin{equation*}
\|u\|_{L^p_s} = \|(I-\Delta)^{s/2} u\|_{L^p} < \infty.
\end{equation*}
Note that $H^s(\mathbb{R}^n) = L^2_s(\mathbb{R}^n)$ with equivalent norms. We denote by $L^2_{\mathrm{comp}}(\mathbb{R}^n)$ and $H^s_{\mathrm{comp}}(\mathbb{R}^n)$ the subspaces of compactly supported elements of $L^2(\mathbb{R}^n)$ and $H^s(\mathbb{R}^n)$, respectively, and by $L^p_{\mathrm{loc}}(\mathbb{R}^n)$, $L^2_{\mathrm{loc}}(\mathbb{R}^n)$, and $H^s_{\mathrm{loc}}(\mathbb{R}^n)$ their local counterparts.
The notation $\|f\|_{L^p_{\mathrm{loc}}}$ is used schematically for the family of seminorms $\|\chi f\|_{L^p}$, where $\chi\in C_0^\infty(\mathbb{R}^n)$. Thus, an estimate $\|f\|_{L^p_{\mathrm{loc}}}\leq C$ means that, for every such $\chi$, there exists a constant $C_\chi$ such that $\|\chi f\|_{L^p}\leq C_\chi$. By $\|T\|_s$ we denote the operator norm of $T$ from $L^2(\mathbb{R}^n)$  to $H^s(\mathbb{R}^n)$.

Let $\mu \in \mathbb{R}$. We recall that $S^{\mu}_{1,0}(\mathbb{R}^n)$ is the space of symbols of order $\mu$ and type $(1,0)$, i.e., $a = a(x,\xi) \in C^{\infty}(\mathbb{R}^n \times \mathbb{R}^n)$ belongs  to $S^{\mu}_{1,0}(\mathbb{R}^n)$ if there exist constants $C_{\alpha,\beta} > 0$ such that
\begin{equation*}
    \forall \alpha,\beta \in \mathbb{N}_0^n : \quad
\left| \partial_x^\alpha \partial_\xi^\beta a(x,\xi)\right|
\leq C_{\alpha, \beta} \langle \xi\rangle^{\mu - |\beta|}
\qquad \forall (x,\xi) \in \mathbb{R}^n \times \mathbb{R}^n.
\end{equation*}
The set of pseudo-differential operators associated to the symbols in $S^{\mu}_{1,0}(\mathbb{R}^n)$ is denoted by $\Psi^{\mu}_{1,0}(\mathbb{R}^n)$. If the symbol has an extra (smooth) dependence on $t \in [0,T]$, we use the notations  $C^{\infty}([0,T], S^\mu)$ and $C^{\infty}([0,T], \Psi^\mu_{1,0})$ for symbols and operators, respectively. For brevity, we adopt the  abbreviated notations $S^\mu$ and $\Psi^\mu$ for  $S^{\mu}_{1,0}(\mathbb{R}^n)$ and $\Psi^{\mu}_{1,0}(\mathbb{R}^n)$,  and $C^{\infty}S^\mu$ and $C^{\infty}\Psi^\mu$ for  $C^{\infty}([0,T], S^\mu)$ and $C^{\infty}([0,T], \Psi^\mu_{1,0})$, respectively.

The space of Fourier integral operators of order $\mu$ with amplitudes of type $(1,0)$ is denoted by $I^\mu$. All Fourier integral operators in this paper are nondegenerate, meaning their canonical relations satisfy the local graph condition.

The key $L^2\to H^s$ estimates for FIOs, the measure estimate for functions with zeros of finite order, and the interpolation lemma are stated and proved in Appendix~\ref{sec:auxiliary}.

\subsection{Microlocal diagonalization}\label{subsec:diag}
We now explain how the microlocal diagonalizability assumption on the principal symbol $A(t,x,\xi)$ of $P$ leads to a global diagonal representation that underlies the Picard series construction.
\begin{remark}\label{rem:covering}
	Since all results of this paper are local in nature, we may assume without loss of generality that the amplitude of $P(t,x,D_x)$ is compactly supported with respect to $x$. Under this assumption, the microlocal diagonalizability of $A(t,x,\xi)$ implies that there exists a finite-sheeted covering on which the eigenvalues $a_j$ and eigenspaces $V_j$ are globally smooth and single-valued; see \cite{KR2007, Roz80}. We henceforth assume that $A$ is smoothly diagonalizable over $\mathbb{R}^n$ and work directly with the diagonal representation.
\end{remark}
Following the reduction described in Remark~\ref{rem:covering}, we write
\begin{equation}\label{eq:decomp}
P(t,x,D_x) = A(t,x,D_x) + B(t,x,D_x),
\end{equation}
where
\begin{equation*}
A(t,x,D_x) = \mathrm{diag}\{A_1(t,x,D_x), \ldots, A_m(t,x,D_x)\}, 
\end{equation*}
is a diagonal matrix of scalar pseudo-differential operators of order one with principal symbols $a_j(t,x,\xi)\in C^{\infty}([0,T]; S^{1}(\mathbb{R}^{n}))$, $j= 1, \ldots, m$. The eigenvalues $a_j$ are not assumed to be distinct: some may be identically equal, or may intersect with finite order of contact governed by Condition~C. We thus allow characteristics of variable multiplicity, including blocks of constant multiplicity as a special case. The operator
\begin{equation*}
B(t,x,D_x) = (B_{jk}(t,x,D_x))_{1 \leq j,k \leq m} 
\in C^\infty([0,T]; \Psi^0(\mathbb{R}^n)^{m \times m})
\end{equation*}
is an $m \times m$ matrix of pseudo-differential operators of order zero. We absorb the diagonal entries of $B$ into the scalar operators by setting
\begin{equation*}
\widehat A_j := A_j+B_{jj}, \qquad
\widehat A := \operatorname{diag}\{\widehat A_1,\ldots,\widehat A_m\}, \qquad
\widehat B := B-\operatorname{diag}\{B_{11},\ldots,B_{mm}\}.
\end{equation*}
Then $P=\widehat A+\widehat B$ and $\widehat B_{jj}=0$ for $1\leq j\leq m$. Since $B_{jj}$ has order zero, $\widehat A_j$ has
the same real principal symbol $a_j$ as $A_j$, so Condition~C is unchanged. The lower-order terms of $\widehat A_j$ include $B_{jj}$; no self-adjointness of $\widehat A_j$ is assumed.

\subsection{The diagonal propagator and the Picard series}
\label{subsec:propagator}

For each $j = 1, \ldots, m$, let $E_j(t, s)$ denote the scalar propagator associated with $\widehat A_j$ solving
\begin{equation}\label{eq:propagator}
\begin{cases}
i\partial_t E_j(t,\tau) = \widehat A_j(t,x,D_x) E_j(t,\tau), \\
E_j(\tau,\tau) = I,
\end{cases}
\end{equation}
for $0 \leq \tau, t \leq T$. By standard theory, $E_j(t,\tau)$ is a nondegenerate Fourier integral operator of order zero. Both $E_j(t,\tau)$ and its inverse are bounded on $L^2(\mathbb{R}^n)$,
uniformly for $t,\tau\in[0,T]$, with $E_j(t,\tau)^{-1}=E_j(\tau,t)$ and evolution property
\begin{equation*}
E_j(t,\tau) E_j(\tau,r) = E_j(t,r), \quad 0 \leq r, \tau, t \leq T.
\end{equation*}
We set $E(t,\tau) = \mathrm{diag}\{E_1(t,\tau),\ldots,E_m(t,\tau)\}$.

Let $U(t,\tau)$ denote the fundamental solution of \eqref{eq:CP}, satisfying
\begin{equation}\label{eq:fundamental}
\begin{cases}
i\partial_t U(t,\tau) = P(t,x,D_x) U(t,\tau), \\
U(\tau,\tau) = I.
\end{cases}
\end{equation}
We factor $U(t,\tau) = E(t,\tau) V(t,\tau)$. Substituting into \eqref{eq:fundamental} and using \eqref{eq:propagator}, we find that $V(t,\tau)$ satisfies
\begin{equation}\label{eq:V}
\begin{cases}
\partial_t V(t,\tau) = Z(t,\tau) V(t,\tau), \\
V(\tau,\tau) = I,
\end{cases}
\end{equation}
where the entries of the matrix-valued operator $Z(t,\tau)$ are given by
\begin{equation}\label{eq:Z}
Z_{jk}(t,\tau) = -i E_j(\tau,t) \widehat B_{jk}(t,x,D_x) E_k(t,\tau),
\quad 1 \leq j,k \leq m.
\end{equation}
Writing the Picard series for \eqref{eq:V}, we obtain
\begin{equation}\label{eq:Picard}
V(t,\tau) = I + \sum_{\ell=1}^{\infty} Q_\ell(t,\tau),
\end{equation}
where the $\ell$-th term is the time-ordered integral
\begin{equation}\label{eq:Ql}
Q_\ell(t,\tau) = \int_\tau^t \int_\tau^{t_1} \cdots \int_\tau^{t_{\ell-1}} 
Z(t_1,\tau) Z(t_2,\tau) \cdots Z(t_\ell,\tau) \, dt_\ell \cdots dt_1.
\end{equation}
The series \eqref{eq:Picard} converges in operator norm on 
$\mathcal{L}(L^2(\mathbb{R}^n))$; indeed, since $Z(t,\tau)$ is bounded on $L^2$ uniformly in $t,\tau \in [0,T]$, the $\ell$-th term satisfies
\begin{equation*}
\|Q_\ell(t,\tau)\|_{L^2 \to L^2} \leq 
\frac{(C(t-\tau))^\ell}{\ell!} \to 0 \quad \text{as } \ell \to \infty.
\end{equation*}
The goal of Section~\ref{sec:picard} is to show that under Condition~C, each operator $Q_\ell$ gains regularity as $\ell$ increases, in the precise sense of Theorem~\ref{thm:smoothing}. In the subsequent analysis, we write $B$ for the off-diagonal remainder $\widehat B$.

\section{Analysis of the Picard Series}\label{sec:picard}
In this section we prove Theorem~\ref{thm:smoothing}, which asserts that the terms $Q_\ell$ of the Picard series \eqref{eq:Picard} gain regularity as $\ell$ increases. The argument proceeds in three stages: separation of  amplitude and phase (Stage~I), construction of the phase function (Stage~II),  and the smoothing estimates themselves (Stage~III).

Recall from \eqref{eq:Ql} that the $\ell$-th Picard iterate is the time-ordered integral
\begin{equation}\label{eq:Ql-recall}
Q_\ell = (-i)^\ell 
\int_0^t \int_0^{t_1} \cdots \int_0^{t_{\ell-1}} 
\mathcal{I}_\ell(t_1,\ldots,t_\ell) \, dt_\ell \cdots dt_1,
\end{equation}
where
\begin{equation*}
	\mathcal{I}_\ell(t_1,\ldots,t_\ell) = E(0,t_1)B(t_1)E(t_1,t_2)\cdots E(t_{\ell-1}, t_{\ell}) B(t_{\ell})E(t_{\ell},0).
\end{equation*}
\subsection{Stage I: Separation of amplitude and phase}
\label{subsec:separation}

The first step is to rearrange each integrand $\mathcal{I}_\ell$ into the  product of a pseudo-differential amplitude operator and a Fourier integral operator carrying the oscillations. The key device is Egorov's theorem (Theorem~\ref{thm:egorov}): for any propagator $E_k$ and any zeroth-order operator $B$,
\begin{equation}\label{eq:conjugation}
E_k(0,t) B(t) 
= \bigl[E_k(0,t)\, B(t)\, E_k(t,0)\bigr] E_k(0,t) 
= \widetilde{B}(t)\, E_k(0,t),
\end{equation}
where $\widetilde{B}(t) = E_k(0,t)\, B(t)\, E_k(t,0) \in C^{\infty}\Psi^0$ is the symbol of $B(t)$ transported backwards from time $t$ to time $0$ along the bicharacteristics of $A_k$.

\subsubsection*{Illustration for $\ell = 1$.}
For a single interaction between modes $j_1$ and $j_2$ at time $t_1$, the integrand reads
\begin{equation*}
\mathcal{I}_1 = E_{j_1}(0,t_1)\, B_{j_1 j_2}(t_1)\, E_{j_2}(t_1,0).
\end{equation*}
Applying \eqref{eq:conjugation} with $k = j_1$, we obtain
\begin{equation*}
\mathcal{I}_1 = \widetilde{B}_1(t_1)\, H(t_1),
\end{equation*}
where
\begin{equation*}
\widetilde{B}_1(t_1) = E_{j_1}(0,t_1)\, B_{j_1 j_2}(t_1)\, E_{j_1}(t_1,0) 
\in C^{\infty}\Psi^0
\end{equation*}
and
\begin{equation*}
H(t_1) = E_{j_1}(0,t_1)\, E_{j_2}(t_1,0)
\end{equation*}
is a Fourier integral operator carrying the oscillatory content.

\subsubsection*{Illustration for $\ell = 2$.}
For two interactions between modes $j_1, j_2, j_3$ at times $t_1 > t_2$, the integrand is
\begin{equation*}
\mathcal{I}_2 = E_{j_1}(0,t_1)\, B_{j_1 j_2}(t_1)\, E_{j_2}(t_1,t_2)\, 
B_{j_2 j_3}(t_2)\, E_{j_3}(t_2,0).
\end{equation*}
Applying \eqref{eq:conjugation} to $B_{j_1 j_2}(t_1)$ as above yields
\begin{equation*}
\mathcal{I}_2 = \widetilde{B}_1(t_1)\, 
\underbrace{E_{j_1}(0,t_1)\, E_{j_2}(t_1,t_2)}_{\Phi}\, 
B_{j_2 j_3}(t_2)\, E_{j_3}(t_2,0).
\end{equation*}
To extract $B_{j_2 j_3}(t_2)$, we conjugate it past the cumulative propagator $\Phi$:
\begin{equation*}
\Phi\, B_{j_2 j_3}(t_2) 
= \bigl[\Phi\, B_{j_2 j_3}(t_2)\, \Phi^{-1}\bigr]\, \Phi 
= \widetilde{B}_2(t_1,t_2)\, \Phi,
\end{equation*}
so that
\begin{equation*}
\mathcal{I}_2 = \widetilde{B}_1(t_1)\, \widetilde{B}_2(t_1,t_2)\, H(t_1,t_2),
\end{equation*}
with
\begin{equation*}
H(t_1,t_2) = E_{j_1}(0,t_1)\, E_{j_2}(t_1,t_2)\, E_{j_3}(t_2,0)
\end{equation*}
and
\begin{equation*}
\widetilde{B}_2(t_1,t_2) 
= E_{j_1}(0,t_1)\, 
\bigl[E_{j_2}(t_1,t_2)\, B_{j_2 j_3}(t_2)\, E_{j_2}(t_2,t_1)\bigr]\, 
E_{j_1}(t_1,0),
\end{equation*}
representing the symbol of $B_{j_2 j_3}$ transported from $t_2$ to $t_1$ along the flow of $a_{j_2}$, and then from $t_1$ to $0$ along the flow of $a_{j_1}$.

\subsubsection*{General formula.}
Iterating this procedure, we arrive at the decomposition
\begin{equation}\label{eq:separation}
\mathcal{I}_\ell(\bar{t}) = \mathcal{B}(\bar{t})\, H(\bar{t}),
\end{equation}
where $\bar{t} = (t_1,\ldots,t_\ell)$ and:
\begin{itemize}
	\item The \emph{phase operator}
	\begin{equation}\label{eq:Hbar}
	H(\bar{t}) = E_{j_1}(0,t_1)\, E_{j_2}(t_1,t_2) \cdots E_{j_{\ell+1}}(t_\ell,0)
	\end{equation}
	is a Fourier integral operator of order zero, with $1 \leq j_k \leq m$ and $j_k \neq j_{k+1}$ for $k = 1,\ldots,\ell$.
	\item The \emph{amplitude operator}
	\begin{equation}\label{eq:Bbar}
	\mathcal{B}(\bar{t}) = \widetilde{B}_1(t_1)\, \widetilde{B}_2(t_1,t_2) 
	\cdots \widetilde{B}_\ell(t_1,\ldots,t_\ell),
	\end{equation}
	with the $k$-th factor given by
	\begin{equation}\label{eq:Btilde-k}
	\widetilde{B}_k = \Biggl(\prod_{r=1}^{k} E_{j_r}(t_{r-1},t_r)\Biggr) 
	B_{j_k j_{k+1}}(t_k) 
	\Biggl(\prod_{r=1}^{k} E_{j_r}(t_{r-1},t_r)\Biggr)^{\!-1},
	\end{equation}
	where $t_0 = 0$, is a pseudo-differential operator of order zero 
	depending smoothly on $\bar{t}$.
\end{itemize}

\subsection{Stage II: Construction of the phase function}
\label{subsec:phase}

Substituting the decomposition \eqref{eq:separation} into \eqref{eq:Ql-recall}, the $\ell$-th Picard iterate takes the form 
\begin{equation}\label{eq:Ql-separated}
Q_\ell = \int_0^t \int_0^{t_1} \cdots \int_0^{t_{\ell-1}}
\mathcal{B}(\bar{t})\, H(\bar{t}) \, dt_\ell \cdots dt_1,
\end{equation}
where $\mathcal{B}(\bar{t}) \in C^{\infty}\Psi^0$ depends smoothly on $\bar{t}$ and $H(\bar{t})$ is given by \eqref{eq:Hbar}. By the theorem on composition of Fourier integral operators, $H(\bar{t})$ is a (parameter-dependent) Fourier integral operator of order zero whose canonical relation $\Lambda^{\bar{t}} \subseteq T^*\mathbb{R}^n \times
T^*\mathbb{R}^n$ is the graph of the composite canonical transformation
$\Psi^{\bar{t}}$:
\begin{equation}\label{eq:canonical-relation}
\Lambda^{\bar{t}} = \bigl\{(x,p,y,\xi) :
(x,p) = \Psi^{\bar{t}}(y,\xi)\bigr\},
\end{equation}
where
\begin{equation}\label{eq:composite-flow}
\Psi^{\bar{t}} = \Phi_{j_1}^{0,t_1} \circ \Phi_{j_2}^{t_1,t_2}
\circ \cdots \circ \Phi_{j_\ell}^{t_{\ell-1},t_\ell}
\circ \Phi_{j_{\ell+1}}^{t_\ell,0},
\end{equation}
and $\Phi_k^{r, \tau}$ denotes the canonical transformation obtained by solving Hamilton's equations
\begin{equation}\label{eq:Hamilton-flow}
\frac{dx}{dt} = \frac{\partial a_k}{\partial \xi}(t,x,\xi), \qquad
\frac{d\xi}{dt} = -\frac{\partial a_k}{\partial x}(t,x,\xi),
\end{equation}
from time $\tau$ to time $r$. When the eigenvalues $a_k$ are independent of $t$, the flow reduces to a one-parameter group $\Phi_k^{r,\tau} = \Phi_k^{r -\tau}$.

\subsubsection*{What we have, what we need.}
The composition theorem for Fourier integral operators identifies $H(\bar t)$ as an FIO of order zero whose canonical relation $\Lambda^{\bar t}$ is the graph of the composite flow $\Psi^{\bar t}$ in~\eqref{eq:composite-flow}. Existence and the geometric structure of the wavefront-set propagation are therefore known abstractly. What is \emph{not} yet known is a concrete oscillatory integral representation with an explicit phase $\phi$ and amplitude $a$, which is what Stage~III requires in order to perform integration by parts in $\xi$ or in $\bar t$ and to estimate the resulting sublevel sets. The role of Stage~II is to construct $\phi$ explicitly. The amplitude $a$ then arises from the standard transport-equation hierarchy associated with the phase, and we keep track only of its principal part throughout.

\subsubsection*{The operator-level identity.}
Differentiating $H(\bar t)=E_{j_1}(0,t_1)\cdots E_{j_{\ell+1}}(t_\ell,0)$ termwise (using $i\partial_t E_j=\widehat A_jE_j$ and the group property $E_j(t,s)E_j(s,r)=E_j(t,r)$) yields
\begin{equation}\label{eq:H-system}
\frac{\partial H(\bar{t})}{\partial t_k}
= i\, T_k(t_1,\ldots,t_k)\, H(\bar{t}),
\quad k = 1,\ldots,\ell,
\end{equation}
where $T_k$ is a pseudo-differential operator of order one given by
\begin{equation}\label{eq:Tk}
T_k(t_1,\ldots,t_k) = E_{j_1}(0,t_1) \cdots E_{j_k}(t_{k-1},t_k)\,
\bigl[\widehat A_{j_k}(t_k) - \widehat A_{j_{k+1}}(t_k)\bigr]\,
E_{j_k}(t_k,t_{k-1}) \cdots E_{j_1}(t_1,0).
\end{equation}
By Egorov's theorem (Lemma~\ref{thm:egorov}), the principal symbol of $T_k$ is
\begin{equation}\label{eq:Tk-symbol}
T_k^0(\bar t,x,\xi)
=
\bigl(a_{j_k}-a_{j_{k+1}}\bigr)\!
\Bigl(t_k,\;
\Phi_{j_k}^{t_k,t_{k-1}}
\circ \Phi_{j_{k-1}}^{t_{k-1},t_{k-2}}
\circ \cdots \circ \Phi_{j_1}^{t_1,0}(x,\xi)
\Bigr),
\quad (x,\xi)\in T^*\mathbb{R}^n,
\end{equation}
which depends only on $(t_1,\ldots,t_k)$ among the time variables.

\subsubsection*{The Hamilton--Jacobi system and the canonical relation.}
We seek a generating function $\psi(\bar t,x,\xi)$ such that $H(\bar t)$ admits the oscillatory representation with phase
\begin{equation*}
\phi(\bar t,x,y,\xi)=\psi(\bar t,x,\xi)-y\cdot\xi.
\end{equation*}
Passing to the principal symbol in \eqref{eq:H-system} via \eqref{eq:Tk-symbol} gives the Hamilton--Jacobi system
\begin{equation}\label{eq:HJ}
\frac{\partial \psi}{\partial t_k}
= T_k^0\Bigl(\bar{t},x,\frac{\partial \psi}{\partial x}\Bigr),
\quad k = 1,\ldots,\ell, \qquad
\psi(0,x,\xi) = x \cdot \xi.
\end{equation}
Its Frobenius integrability is established in Lemma~\ref{lem:frobenius} below. The canonical relation $\Lambda^{\bar t}$ is parametrized via
\begin{equation}\label{eq:Lagrangian}
y = \frac{\partial \psi}{\partial \xi}(\bar t,x,\xi),\qquad
p = \frac{\partial \psi}{\partial x}(\bar t,x,\xi),
\end{equation}
and we write $(x^{\bar t}(y,\xi),\,p^{\bar t}(y,\xi))=\Psi^{\bar t}(y,\xi)$ for the resulting image.

\begin{lemma}[Motion of the canonical-relation endpoint]\label{lem:hamilton-eqs}
For each $k=1,\ldots,\ell$, the image point $(x^{\bar t},p^{\bar t})=\Psi^{\bar t}(y,\xi)$ satisfies
\begin{equation}\label{eq:Hamilton-eqs}
\frac{\partial x^{\bar{t}}}{\partial t_k}
= -\frac{\partial T_k^0}{\partial p}\bigl(\bar t,x^{\bar t},p^{\bar t}\bigr), \qquad
\frac{\partial p^{\bar{t}}}{\partial t_k}
= +\frac{\partial T_k^0}{\partial x}\bigl(\bar t,x^{\bar t},p^{\bar t}\bigr).
\end{equation}
\end{lemma}

\begin{proof}
Differentiate $y=\partial_\xi\psi(\bar t,x^{\bar t},\xi)$ in $t_k$ at fixed $(y,\xi)$:
\[
0=\partial_{t_k}\partial_\xi\psi\bigr|_{(x,\xi)}+\partial_x\partial_\xi\psi\cdot\partial_{t_k}x^{\bar t}.
\]
By \eqref{eq:HJ}, $\partial_{t_k}\partial_\xi\psi\bigr|_{(x,\xi)}=\partial_\xi\bigl[T_k^0(\bar t,x,\partial_x\psi)\bigr]=(T_k^0)_p\,\partial_x\partial_\xi\psi$. Since $\det\partial_x\partial_\xi\psi\ne 0$ for small $\bar t$, this yields $\partial_{t_k}x^{\bar t}=-(T_k^0)_p$. Differentiating $p=\partial_x\psi$ in $t_k$ and using $\partial_x[T_k^0(\bar t,x,\partial_x\psi)]=(T_k^0)_x+(T_k^0)_p\,\partial_x^2\psi$ together with the just-proved expression for $\partial_{t_k}x^{\bar t}$, the cross terms cancel and $\partial_{t_k}p^{\bar t}=+(T_k^0)_x$.
\end{proof}

\begin{remark}
The minus sign in the first equation of \eqref{eq:Hamilton-eqs} reflects the fact that $\Lambda^{\bar t}$ is the canonical relation of an FIO with phase $e^{i(\psi-y\cdot\xi)}$: the image point flows along the bicharacteristic of $T_k^0$ in the \emph{reversed} time direction. As a consistency check, take $\ell=1$, $a_{j_1} \equiv 0$, $a_{j_2 }= a$. Then $H(t_1)=E_a(t_1,0)$ and $\Psi^{t_1}=\Phi_a^{t_1,0}$ (forward $a$-flow), while $T_1^0=-a$; equation~\eqref{eq:Hamilton-eqs} then gives $\partial_{t_1}x=-(T_1^0)_p=+a_p$, in agreement with the bicharacteristic of $a$.
\end{remark}

\begin{lemma}[Time-dependent Frobenius identity]\label{lem:frobenius}
For $ 1 \leq k < j \leq \ell$,
\begin{equation}\label{eq:Frobenius}
\frac{\partial T_j^0}{\partial t_k}(\bar t, x, \xi) \;=\;\{T_k^0,T_j^0\}(\bar t,x,\xi).
\end{equation}
In particular, the system \eqref{eq:HJ} is integrable for small $\bar t$, and the resulting phase function $\psi$ is non-degenerate.
\end{lemma}

\begin{proof}
The equality of mixed partials $\partial_{t_j}\partial_{t_k}H=\partial_{t_k}\partial_{t_j}H$ applied to \eqref{eq:H-system} yields
\begin{equation}\label{eq:compatibility-op}
i\Bigl(\frac{\partial T_j}{\partial t_k}-\frac{\partial T_k}{\partial t_j}\Bigr)\;=\;[T_j,T_k]\quad\bmod\Psi^0.
\end{equation}
Taking principal symbols and using $[T_j,T_k]^0=-i\{T_j^0,T_k^0\}$,
\[
\frac{\partial T_j^0}{\partial t_k}-\frac{\partial T_k^0}{\partial t_j}\;=\;-\{T_j^0,T_k^0\}\;=\;\{T_k^0,T_j^0\}.
\]
For $j>k$, equation~\eqref{eq:Tk-symbol} shows that $T_k^0$ depends only on $(t_1,\ldots,t_k)$, hence $\partial T_k^0/\partial t_j=0$ and \eqref{eq:Frobenius} follows. Integrability of \eqref{eq:HJ} is the same statement at the level of the symbol equation, so $\psi$ exists and is non-degenerate for $\bar t$ in a neighbourhood of $0$.
\end{proof}

\begin{remark}
In the autonomous case $\partial_t a_j\equiv 0$, the symbol $T_k^0$ depends on $\bar t$ only through the flow durations $t_k-t_{k-1}$, and \eqref{eq:Frobenius} reduces to the identity established by Rozenblum \cite{Roz80}. The novelty in the present setting is that the explicit time-dependence of $a_j(t,x,\xi)$ contributes new terms inside $T_k^0$; the Frobenius condition \eqref{eq:Frobenius} nevertheless persists unchanged because it is a structural consequence of the operator equation~\eqref{eq:H-system}, independent of whether the eigenvalues depend on $t$ or not.
\end{remark}

\subsubsection*{Derivatives of the symbol along the flow.}
For fixed $(y,\xi)$, define
\begin{equation}\label{eq:f-def}
f_j(\bar{t}) = T_j^0\bigl(\bar{t},\, x^{\bar{t}}(y,\xi),\, p^{\bar{t}}(y,\xi)\bigr),\qquad j=1,\ldots,\ell.
\end{equation}
For each $k$, denote by
\begin{equation}\label{eq:zk-def}
z_k(\bar t,x,\xi)\;=\;\Phi_{j_k}^{t_k,t_{k-1}}\!\circ\Phi_{j_{k-1}}^{t_{k-1},t_{k-2}}\!\circ\cdots\circ\Phi_{j_1}^{t_1,0}(x,\xi)
\end{equation}
the partial composite endpoint appearing in \eqref{eq:Tk-symbol}, so that
$T_k^0(\bar t,x,\xi)=(a_{j_k}-a_{j_{k+1}})(t_k,z_k(\bar t,x,\xi))$.

\begin{proposition}\label{prop:fj-derivative}
The function $f_j(\bar t)$ defined in \eqref{eq:f-def} satisfies:
\begin{enumerate}
\item[(i)] (\emph{Invariance.}) For $1\le k<j\le\ell$, $\displaystyle \partial_{t_k} f_j \equiv 0.$
\item[(ii)] (\emph{Self-derivative.}) For $1\le k\le\ell$,
\[
\partial_{t_k} f_k(\bar t)
\;=\;\Bigl[\,\partial_t(a_{j_k}-a_{j_{k+1}})\;+\;\{a_{j_{k+1}},a_{j_k}\}\,\Bigr]\bigl(t_k,\,z_k(\bar t,x^{\bar t},p^{\bar t})\bigr).
\]
\item[(iii)] (\emph{Iteration.}) For every $\nu\ge 1$, $\partial_{t_k}^{\nu}f_k=\mathcal{D}_{j_k}^{\nu}(a_{j_k}-a_{j_{k+1}})$ evaluated at $(t_k,z_k(\bar t,x^{\bar t},p^{\bar t}))$.
\end{enumerate}
\end{proposition}

\begin{proof}
\emph{(i) Case $k<j$.}\;
By the chain rule together with Lemmas~\ref{lem:hamilton-eqs} and~\ref{lem:frobenius},
\begin{align*}
\frac{\partial f_j}{\partial t_k}
&=\;\frac{\partial T_j^0}{\partial t_k}\Big|_{(x,p)}+\frac{\partial T_j^0}{\partial x}\frac{\partial x^{\bar t}}{\partial t_k}+\frac{\partial T_j^0}{\partial p}\frac{\partial p^{\bar t}}{\partial t_k}\\
&=\;\{T_k^0,T_j^0\}\;+\;\frac{\partial T_j^0}{\partial x}\!\cdot\!\Bigl(-\tfrac{\partial T_k^0}{\partial p}\Bigr)\;+\;\frac{\partial T_j^0}{\partial p}\!\cdot\!\Bigl(\tfrac{\partial T_k^0}{\partial x}\Bigr)\\
&=\;\{T_k^0,T_j^0\}\;-\;\Bigl(\tfrac{\partial T_j^0}{\partial x}\tfrac{\partial T_k^0}{\partial p}-\tfrac{\partial T_j^0}{\partial p}\tfrac{\partial T_k^0}{\partial x}\Bigr)\;=\;\{T_k^0,T_j^0\}-\{T_k^0,T_j^0\}\;=\;0.
\end{align*}

\emph{(ii) Case $k=j$.}\;
Since $\{T_k^0,T_k^0\}=0$, the implicit terms vanish and
\begin{equation*}
    \frac{\partial f_k}{\partial t_k}\;=\;\frac{\partial T_k^0}{\partial t_k}\Big|_{(x,p)}.
\end{equation*}
The right-hand side is a derivative \emph{at fixed $(x,p)$} of $T_k^0(\bar t,x,p)=(a_{j_k}-a_{j_{k+1}})(t_k,z_k)$. The dependence on $t_k$ enters in two ways: explicitly through the time slot of $a_{j_k}-a_{j_{k+1}}$, and implicitly through the partial endpoint $z_k$, whose evolution with $t_k$ is governed by the \emph{outermost} flow segment $\Phi_{j_k}^{t_k,t_{k-1}}$ alone:
\begin{equation*}
    \partial_{t_k}z_k\;=\;\bigl((a_{j_k})_p,\,-(a_{j_k})_x\bigr)\bigr|_{(t_k,z_k)}.
\end{equation*}
Hence
\begin{equation*}
    \frac{\partial T_k^0}{\partial t_k}\Big|_{(x,p)}
\;=\;\partial_t(a_{j_k}-a_{j_{k+1}})\;+\;\{a_{j_k},\,a_{j_k}-a_{j_{k+1}}\}
\;=\;\partial_t(a_{j_k}-a_{j_{k+1}})\;+\;\{a_{j_{k+1}},\,a_{j_k}\},
\end{equation*}
both evaluated at $(t_k,z_k)$, where we used $\{a_{j_k},a_{j_k}\}=0$ and $\{a_{j_k},-a_{j_{k+1}}\}=\{a_{j_{k+1}},a_{j_k}\}$. This is~(ii).

\emph{(iii) Iteration.}\;
Each subsequent application of $\partial_{t_k}$ is a derivative of a function $g(t_k,z_k(t_k))$ of the same form as $T_k^0$. The chain rule used in~(ii) gives $\partial_{t_k}g=\partial_t g+\{a_{j_k},g\}$ at $(t_k,z_k)$. Applying this $\nu-1$ times to the result of (ii) proves (iii).
\end{proof}

\begin{remark}\label{rem:condC-bridge}
Comparing Proposition~\ref{prop:fj-derivative}\,(iii) with the definition of Condition~C in~\eqref{eq:ConditionC},
\begin{equation*}
    \partial_{t_k}^{\nu}f_k(\bar t)\;=\;\mathcal{D}_{j_k}^{\nu}(a_{j_k}-a_{j_{k+1}})\bigl(t_k,\,z_k(\bar t,x^{\bar t},p^{\bar t})\bigr).
\end{equation*}
Consequently, under Condition~C every zero of the function $t_k\mapsto f_k(\bar t)$ has order at most $M$, uniformly in the remaining variables. This is precisely the input to Proposition~\ref{prop:measure} that yields the measure estimate~\eqref{eq:meas_estimate} in Stage~III.
\end{remark}

\subsection{Stage III: Smoothing estimates}
\label{subsec:smoothing}

We now establish the smoothing properties of $Q_\ell$. From the preceding analysis, $Q_\ell$ can be represented as an oscillatory integral operator
\begin{equation}\label{eq:Ql-osc}
Q_\ell u(x) = \int_\Delta \int_{\mathbb{R}^n} \int_{\mathbb{R}^n} 
e^{i\phi(\bar{t},x,y,\xi)}\, b(\bar{t},x,y,\xi)\, u(y)\, 
dy\, d\xi\, d\bar{t},
\end{equation}
where $\phi$ is the phase function satisfying \eqref{eq:HJ} and 
\eqref{eq:Lagrangian}, $b$ is an amplitude of order zero (compactly  supported in $x$ and $y$), and $\Delta = \{0 \leq t_\ell \leq \cdots \leq t_1 \leq t\}$ is the standard simplex in $[0,t]^{\ell}$.

\subsubsection*{Three-region decomposition.}
Let $\chi\in C_0^{\infty}(\mathbb{R})$ be a cutoff with $\chi(\tau)=1$ for $|\tau|<1$ and $\chi(\tau)=0$ for $|\tau|>2$.
Fix $\varepsilon>0$ small. We decompose
\begin{equation*}
    Q_{\ell} \;=\; R_1 + R_2 + R_3,
\end{equation*}
where
\begin{equation*}
R_j u(x) = \int_\Delta \int_{\mathbb{R}^n} \int_{\mathbb{R}^n} 
e^{i\phi(\bar{t},x,y,\xi)}\, \mu_j(\varepsilon,\bar{t},x,y,\xi)\, 
u(y)\, dy\, d\xi\, d\bar{t},
\end{equation*}
with cutoff functions
\begin{align}
\mu_1 &= \Bigl(1 - \chi\bigl(\varepsilon^{-1}|\partial_\xi\phi|\bigr)\Bigr)\, b(\bar{t},x,y,\xi), 
\label{eq:mu1} \\
\mu_2 &= \chi\bigl(\varepsilon^{-1}|\partial_\xi\phi|\bigr)\, 
\Bigl(1 - \chi\bigl(\varepsilon^{-1}|\partial_{\bar{t}}\phi|\,|\xi|^{-1}\bigr)\Bigr)\, b(\bar{t},x,y,\xi), 
\label{eq:mu2} \\
\mu_3 &= \chi\bigl(\varepsilon^{-1}|\partial_\xi\phi|\bigr)\, 
\chi\bigl(\varepsilon^{-1}|\partial_{\bar{t}}\phi|\,|\xi|^{-1}\bigr)\, b(\bar{t},x,y,\xi). 
\label{eq:mu3}
\end{align}
The three regions are characterized as follows:
\begin{align*}
\operatorname{supp} \mu_1 &\colon |\partial_\xi \phi| \gtrsim \varepsilon; \\
\operatorname{supp} \mu_2 &\colon |\partial_\xi \phi| \lesssim \varepsilon 
\text{ and } |\partial_{\bar{t}} \phi| \gtrsim \varepsilon|\xi|; \\
\operatorname{supp} \mu_3 &\colon |\partial_\xi \phi| \lesssim \varepsilon 
\text{ and } |\partial_{\bar{t}} \phi| \lesssim \varepsilon|\xi|.
\end{align*}

\subsubsection*{The Non-Stationary Region $R_1$.} On the support of $\mu_1$, the phase function is non-stationary with respect to $\xi$ with lower-bound
\begin{equation*}
\left| \partial_{\xi} \phi(\bar{t},x,y,\xi) \right| \geq \varepsilon.
\end{equation*}
Therefore, we can construct a first-order differential operator $L(\partial_{\xi})$ satisfying $L e^{i\phi} = e^{i\phi}$, whose coefficients are estimated by $C\varepsilon^{-1}$ on the support of the amplitude of $\mu_1$. We proceed by integrating by parts $p$ times with the operator $L$. In each step, differentiation introduces a factor $\varepsilon^{-1}$ from the coefficients of $L$, and an additional factor of $\varepsilon^{-1}$ arises from differentiating the scaled cutoff function $\mu_1$ (where $\partial_{\xi} \mu_1 \sim \varepsilon^{-1}$). Thus, after $p$ iterations, we obtain an operator with an amplitude of order $-p$ and coefficients bounded by $\varepsilon^{-2p}$. Applying Lemma~\ref{lem:FIO-est} with $q=-p$, we arrive at the following estimate 
\begin{equation}\label{eq:R1-est}
\| R_1 \|_{p-n-1} \leq C \varepsilon^{-3p + n+1}.
\end{equation}

\subsubsection*{The Time-Non-Stationary Region $R_2$.} On the support of $ \mu_2$, the phase function may be stationary in $\xi$, but it is non-stationary with respect to the time variables: 
\begin{equation*}
| \partial_{\xi} \phi(\bar{t},x,y,\xi) | \leq 2\varepsilon, \quad |\partial_{\bar t} \phi(\bar{t},x,y,\xi) | \geq \varepsilon |\xi|.
\end{equation*}
Therefore, the procedure of integrating by parts with respect to $\xi$ cannot be applied to $R_2$. However, we can instead apply integration by parts with respect to $\bar{t}$. We construct a differential operator $M(\partial_{\bar{t}}) $ satisfying $M e^{i\phi} = e^{i\phi}$. On the support of $\mu_2$, the coefficients of $M$ are bounded by $C \varepsilon^{-1}|\xi|^{-1}$ (the factor $|\xi|^{-1}$ arises because $\partial_{\bar{t}} \phi$ scales linearly with $\xi$). Performing integration by parts with $M$ yields the following trade-off:
\begin{itemize}
	\item \textit{Amplitude Decay:} The order of the amplitude decreases by $1$ (due to the $|\xi|^{-1}$ factor).
	\item \textit{Coefficient Growth:} The bound increases by a factor of $\varepsilon^{-2}$ (one $\epsilon^{-1}$ from the operator $M$ and another $\varepsilon^{-1}$ from differentiating the cutoff function).
\end{itemize}
The boundary integrals have the same form as $Q$ but they have amplitudes of order $-1$ and depend on not less than $\ell -2$ time variables.  The reason for possibly losing two variables is that after restriction to the boundary, say $t_2=t_1$, it may happen that $a_{j_1}$ and $a_{j_3}$ are the same roots. It follows that this procedure can be iterated at most $[\ell/2]$ times.  After $[\ell /2]$ steps, we obtain operators of order $-[\ell/2]$ with coefficients bounded by $C \varepsilon^{-2[\ell/2]}$. Applying Lemma~\ref{lem:FIO-est} for $[\ell/2]>n+1$, we obtain the estimate  
\begin{equation}\label{eq:R2-est}
\| R_2 \|_{[\ell/2] - n-1} \leq C \varepsilon^{-3[\ell/2] +n+1}.
\end{equation}

\subsubsection*{The Critical Region $R_3$.}
On the support of $\mu_3$, the phase function is nearly stationary in all variables, so integration by parts in $\xi$ or in $\bar t$ is unavailable; in particular,
\begin{equation}\label{eq:zones}
|\partial_\xi\phi(\bar t,x,y,\xi)|\leq 2\varepsilon,\qquad |\partial_{\bar t}\phi(\bar t,x,y,\xi)|\leq 2\varepsilon|\xi|.
\end{equation}
Since the singular support of the integral kernel of $R_3$ may be highly irregular here, a more delicate analysis is required. The key step is to show that, on the support of $\mu_3$,
\begin{equation}\label{eq:T0_estim}
|T_j^0(\bar t,x^{\bar t}(y,\xi),p^{\bar t}(y,\xi))|\;\leq\; C\,\varepsilon|\xi|,\qquad j=1,\ldots,\ell.
\end{equation}
By \eqref{eq:HJ} and \eqref{eq:zones},
\begin{equation}\label{est_T0}
|T_j^0(\bar t,x,\partial_x\psi)|\;\leq\; C\,\varepsilon|\xi|,\qquad j=1,\ldots,\ell.
\end{equation}
Since $T_j^0$ is positively homogeneous of degree one in $\xi$, we also have
\begin{equation*}
|\partial_x T_j^0(\bar t,x,\partial_x\psi)|\leq C|\xi|, \qquad
|\partial_\xi T_j^0(\bar t,x,\partial_x\psi)|\leq C,\qquad j=1,\ldots,\ell.
\end{equation*}
For $(x^{\bar t}(y,\xi),p^{\bar t}(y,\xi))=\Psi^{\bar t}(y,\xi)$ the canonical-relation endpoint, the mean value theorem and the homogeneity bounds above give
\begin{equation*}
    \bigl|T_j^0(\bar t,x,\partial_x\psi)-T_j^0(\bar t,x^{\bar t}(y,\xi),p^{\bar t}(y,\xi))\bigr|
\;\leq\; C\bigl(|\xi|\,|x-x^{\bar t}(y,\xi)|+|\partial_x\psi-p^{\bar t}(y,\xi)|\bigr).
\end{equation*}
By construction of the phase function and of $\Lambda^{\bar t}$, we have $\partial_\xi\phi(\bar t,x^{\bar t}(y,\xi),y,\xi)=0$. Writing the difference as an integral remainder,
\begin{equation*}
    \partial_\xi\phi(\bar t,x,y,\xi)
\;=\;\Bigl(\int_0^1\partial^2_{x,\xi}\phi\bigl(\bar t,x^{\bar t}+s(x-x^{\bar t}),y,\xi\bigr)\,ds\Bigr)(x-x^{\bar t}),
\end{equation*}
and using $\det\partial^2_{x,\xi}\phi\ne 0$ for small $\bar t$, we obtain
\begin{equation*}
    |x-x^{\bar t}(y,\xi)|\;\leq\; C\,|\partial_\xi\phi(\bar t,x,y,\xi)|.
\end{equation*}
Similarly, from the parametrization \eqref{eq:Lagrangian} of $\Lambda^{\bar t}$,
\begin{align}\label{eq:p-est}
\bigl|\partial_x\psi(\bar t,x,\xi)-p^{\bar t}(y,\xi)\bigr|
&= \bigl|p^{\bar t}(\partial_\xi\psi(\bar t,x,\xi),\xi)-p^{\bar t}(y,\xi)\bigr| \nonumber\\
&\leq C\,|\partial_\xi\phi(\bar t,x,y,\xi)|\,|\xi|.
\end{align}
Combining \eqref{eq:zones}--\eqref{eq:p-est} yields the claimed estimate~\eqref{eq:T0_estim} on the support of~$\mu_3$.

Estimate~\eqref{eq:T0_estim} restricts the support of $\mu_3$ to
\begin{equation*}
\Xi_1 \;=\; \Bigl\{(\bar t,x,y,\xi):x\in\mathbb{R}^n,\;\xi\neq 0,\;(\bar t,y,\xi/|\xi|)\in\Xi\Bigr\},
\end{equation*}
where
\begin{equation*}
\Xi \;=\; \Bigl\{(\bar t,y,\eta)\in[0,T]^\ell\times S^*\mathbb{R}^n:\bigl|T_j^0\bigl(\bar t,x^{\bar t}(y,\eta),p^{\bar t}(y,\eta)\bigr)\bigr|\leq C\varepsilon,\;j=1,\ldots,\ell\Bigr\}.
\end{equation*}
We may therefore apply estimate~\eqref{eq:FIO-meas} of Lemma~\ref{lem:FIO-est} with $q=0$ and $\delta=1$ to the operator $R_3$:
\begin{equation}\label{eq:R3-meas}
\|R_3\|_{-3n/2-2}\;\leq\;C\,\varepsilon^{-n-1}\,\mathrm{meas}(\Xi),
\end{equation}
where $\mathrm{meas}$ is the canonical induced measure on $[0,T]^\ell\times S^*\mathbb{R}^n$.

\subsubsection*{The measure estimate via Condition C}
It remains to bound $\mathrm{meas}(\Xi)$. Since the Hamiltonian flows $\Phi_k^{\tau,s}$ preserve the symplectic measure, we may replace $\Xi$ by
\begin{equation*}
\Xi' \;=\; \bigl\{(\bar t,y,\xi):|T_j^0(\bar t,y,\xi)|\leq C\varepsilon,\;c\leq|\xi|\leq C,\;j=1,\ldots,\ell\bigr\}.
\end{equation*}
For each $j$, define the slice
\begin{equation*}
\Sigma_j(t_1,\ldots,\widehat{t_j},\ldots,t_\ell,y,\xi)
\;=\;\bigl\{t_j\in[0,T]:|T_j^0(\bar t,y,\xi)|\leq C\varepsilon\bigr\}.
\end{equation*}
By Remark~\ref{rem:condC-bridge}, the function $t\mapsto T_j^0(t_1,\ldots,t,\ldots,t_\ell,y,\xi)$ has $\partial_t^\nu = \mathcal{D}_{j_k}^\nu(a_{j_k}-a_{j_{k+1}})$ evaluated at the partial composite endpoint, so Condition~C ensures that this function vanishes to order at most~$M$ at each zero, uniformly in the remaining variables (with $(y,\xi)$ in a compact set and $c\leq|\xi|\leq C$). Proposition~\ref{prop:measure} therefore gives
\begin{equation}\label{eq:meas_estimate}
\sup_{w}\,\mathrm{meas}\,\Sigma_j(w,y,\xi)\;\leq\;C\,\varepsilon^{1/2M}.
\end{equation}
Iterating Fubini in the $\ell$ time variables, with~\eqref{eq:meas_estimate} applied at each step, yields
\begin{equation*}
\mathrm{meas}(\Xi)\;\leq\;C\prod_{j=1}^{\ell}\sup_w\mathrm{meas}\,\Sigma_j(w,y,\xi)\;\leq\;C\,\varepsilon^{\ell/2M}.
\end{equation*}
Combining with~\eqref{eq:R3-meas} gives
\begin{align}\label{eq:R3}
\|R_3\|_{-3n/2-2}\;\leq\;C\,\varepsilon^{\ell/2M-n-1}.
\end{align} 

\subsubsection*{Conclusion via Interpolation}

We now combine the three estimates \eqref{eq:R1-est}, \eqref{eq:R2-est}, and \eqref{eq:R3} using Lemma~\ref{lem:interpolation}. For fixed \(\ell\), we choose
\[
p = -\frac{3n}{2}-2,\qquad
r = [\ell/2] - n - 1,
\]
and set
\[
T = 3[\ell/2] - 1 - n, \qquad
S = \frac{\ell}{2M} - n - 1.
\]
For sufficiently large \(\ell\) we have \(S>0\), and estimates
\eqref{eq:R1-est}, \eqref{eq:R2-est}, \eqref{eq:R3} show that for each small \(\varepsilon>0\),
\begin{equation*}
    Q_\ell u = u_1^{(\varepsilon)} + u_2^{(\varepsilon)},
\end{equation*}
with
\begin{equation*}
    \|u_1^{(\varepsilon)}\|_{H^r} \le C\varepsilon^{-T},\qquad
\|u_2^{(\varepsilon)}\|_{H^p} \le C\varepsilon^{S}.
\end{equation*}
By Lemma~\ref{lem:interpolation}, it follows that \(Q_\ell\) is bounded from \(L^2_{\mathrm{comp}}\) to \(H_{\mathrm{loc}}^{N(\ell)}\) for any
\begin{equation*}
    N(\ell) < \frac{pT + rS}{T+S},
\end{equation*}
with
\begin{equation*}
    N(\ell)
= \max \left\{0, \frac{\big(-\tfrac{3n}{2}-2\big)\big(3[\ell/2]-1-n\big) + \big([\ell/2]-n-1\big)\big(\tfrac{\ell}{2M}-n-1\big)}
{3[\ell/2]-2n-2 + \tfrac{\ell}{2M}}\right\}.
\end{equation*}
A straightforward asymptotic analysis as \(\ell\to\infty\) shows that
\begin{equation}\label{eq:Nell}
N(\ell) \sim \frac{\ell}{6M+2}, 
\end{equation}
so in particular \(N(\ell)\to\infty\) with \(\ell \rightarrow \infty\). This proves the smoothing property in Theorem~\ref{thm:smoothing}.
\qed
\section{From Picard Smoothing to Regularity}\label{sec:principle}
The purpose of this section is to pass from the smoothing estimates for the Picard terms to \(L^p\)-regularity of solutions. We first formulate an abstract regularity principle for nonhomogeneous first-order Cauchy problems. We then record two model consequences for pseudo-differential and Fourier integral operator families, and finally apply the principle to first-order hyperbolic systems and to scalar weakly hyperbolic equations.

\subsection{Abstract regularity principle}
We consider the Cauchy problem
\begin{align}\label{eq:3.1}
\begin{cases}
u' - Z(t)u = r(t), \\
u(0) = u_0,
\end{cases}
\end{align}
where the operator \(Z(t)\) is a time-dependent family satisfying
\[
Z(t) \in \mathcal{L}\bigl(C_0^{\infty}(\mathbb{R}^n), \mathcal{D}'(\mathbb{R}^n)\bigr),
\qquad t\in[0,T].
\]
Let \(W_0\), \(W_1\), and \(W\) be linear topological subspaces of \(\mathcal{D}'(\mathbb{R}^n)\) with
\[
W_0,\; W_1 \hookrightarrow W.
\]
We shall assume that the forcing term satisfies \(r(t)\in W_0\) and that the initial datum satisfies \(u_0\in W_1\). The particular choice of these spaces will be made later, according to the properties of \(Z(t)\).

The question is: under which assumptions on \(Z(t)\) does the solution of \eqref{eq:3.1} satisfy
\begin{equation*}
    u(t,\cdot)\in W,\qquad 0\le t\le T,
\end{equation*}
whenever \(r(t)\in W_0\) and \(u_0\in W_1\)?

In applications, the evolution generated by \(Z(t)\) may involve a loss of regularity. The role of the following principle is to isolate the amount of regularity that can be recovered from two sources: the smoothing produced by iterated time integration and the regularity theory for the nonhomogeneous problem with zero initial data.

Since our main concern is regularity rather than existence and uniqueness, we assume the following \(L^2\)-well-posedness property.

\noindent\textbf{(WP)}\quad
For every admissible right-hand side, the Cauchy problem \eqref{eq:3.1} with zero initial datum,
\begin{equation*}
    u(0)=0,
\end{equation*}
is uniquely solvable in the class of functions \(u=u(t,x)\) satisfying
\begin{equation*}
 u(t,\cdot) \in L^2(\mathbb{R}^n),\qquad 0\leq t\leq T.   
\end{equation*}

The Picard representation then allows one to reduce the case of general initial data to the zero-initial-data problem. Thus the assumption \textbf{(WP)} is sufficient for the regularity argument below.

\begin{theorem}\label{thm:3.1}
	Let $W_0, W_1\hookrightarrow W$ be linear subspaces of $\mathcal{D}'(\mathbb{R}^n)$, and let
	\[
	Z(t)\in\mathcal{L}\big(C_0^{\infty}(\mathbb{R}^n),\mathcal{D}'(\mathbb{R}^n)\big),\qquad t\in[0,T].
	\]
	Assume the following.
	\begin{itemize}
		\item[(i)] \textbf{Boundedness.} The family $Z$ extends to an element of
		$L^{\infty}\big([0,T],\mathcal{L}(L^2(\mathbb{R}^n),L^2(\mathbb{R}^n))\big)$.
		Moreover, for each $t\in[0,T]$, the operator $Z(t)$ extends to a continuous linear map $W_1\to W$.
		
		\item[(ii)] \textbf{Stability under composition.} For every $m\in\mathbb{N}$ and every $t_1,\dots,t_m\in[0,T]$, the composition
		$Z(t_1)\circ\cdots\circ Z(t_m)$ extends to a continuous linear map $W_1\to W$.
		
		\item[(iii)] \textbf{Smoothing by time integration.} There exists an integer $\ell\ge 1$ such that, for every $t\in[0,T]$, the operator
		\[
		Z(t) \int_0^{t} \int_0^{t_1} \cdots \int_0^{t_{\ell -1 } } Z(t_1) \circ \cdots \circ Z(t_{\ell})\, dt_{\ell} \cdots dt_{1}
		\]
		extends to a continuous linear map $W_1\to W_0$.
		
		\item[(iv)] \textbf{Regularity for zero Cauchy data.} For each $r(t)\in W_0$, the solution $v$ of
		\begin{align}\label{eq:3.2}
		\begin{cases}
		v' - Z(t)v = r(t), & t\in[0,T],\\
		v(0)=0,
		\end{cases}
		\end{align}
		satisfies (WP) and belongs to $W$ at each time, i.e.
		$v(t,\cdot)\in W$ for all $t\in[0,T]$.
	\end{itemize}
	
	Then, for each $r(t)\in W_0$ and each initial datum $u(0)\in W_1$, the solution $u$ of
	\begin{align}\label{eq:3.3}
	\begin{cases}
	u' - Z(t)u = r(t), & t\in[0,T],\\
	u(0)\in W_1,
	\end{cases}
	\end{align}
	belongs to $W$ at each time: $u(t,\cdot)\in W$ for all $t\in[0,T]$.
	
	If, in addition, $W_0$, $W_1$, and $W$ are normed spaces and the solutions $v$ in \eqref{eq:3.2} satisfy the estimate
	\[
	\|v(t,\cdot)\|_W \le C\|r(t)\|_{W_0},\qquad t\in[0,T],
	\]
	then $u$ satisfies
	\[
	\|u(t,\cdot)\|_W \le C\big(\|u(0)\|_{W_1}+\|r(t)\|_{W_0}\big),\qquad t\in[0,T].
	\]
\end{theorem}
\begin{remark}
	The proof of Theorem~\ref{thm:3.1} was given in \cite{KR2007}. In the present
	paper we use this abstract regularity principle as a black-box tool, applying
	it with choices of \(W_0\), \(W_1\), and \(W\) adapted to the operator family
	under consideration.
\end{remark}

\begin{remark}
	The assumptions in Theorem~\ref{thm:3.1} have separate roles. Conditions
	\emph{(i)} and \emph{(ii)} impose the required mapping stability of the
	operator family \(Z(t)\). In particular, when \(W_1\) is strictly smaller than
	\(W\), condition \emph{(ii)} does not follow from \emph{(i)} alone. In the
	applications below, \(Z(t)\) is a family of time-dependent pseudo-differential
	or Fourier integral operators, and finite compositions of such operators
	remain within the same operator-theoretic framework.
	
	Condition \emph{(iii)} is the smoothing input. It reflects the familiar
	phenomenon from harmonic analysis that integration with respect to a parameter
	can improve regularity. Finally, condition \emph{(iv)} supplies the
	regularity theory for the nonhomogeneous problem with zero Cauchy data.
	Together, these assumptions allow one to transfer the smoothing gained from
	iterated time integration to the solution of the full Cauchy problem.
\end{remark}

\subsection{Proof of Theorem~\ref{thm:Lp-reg}}

We rewrite the system in a form to which Theorem~\ref{thm:3.1} applies.
	
\smallskip
\noindent\textit{Step 1: Reduction to a block-diagonal principal part.}
	By Remark~\ref{rem:covering}, the localization to compactly $x$-supported amplitudes makes the characteristic roots globally smooth and single-valued, and we may work directly with the diagonal representation of $A$. This local lifting argument is valid in every dimension $n\ge 1$ and is independent of the global manifold version (Lemma~\ref{lem:diag}), which would require $n\ge 3$. We thus write
	\begin{align*}
	P(t,x,D_x) &= \bigoplus_{j=1}^m a_j(t,x,D_x) + B(t,x,D_x),\\
	B(t,x,D_x) &= \big(B_{jk}(t,x,D_x)\big)_{1\le j,k\le m},\qquad B_{jk}\in C^{\infty}([0,T],\Psi^0).
	\end{align*}
	
	Some of the roots $a_j$ may coincide identically. We group equal roots into blocks: there exist indices
	\[
	1=k_1<k_2<\cdots<k_r=m+1
	\]
	so that $a_{k_i}\equiv\cdots\equiv a_{k_{i+1}-1}$ for each $i$, and roots from different blocks are not identically equal. This yields a block-diagonal decomposition of the first-order part,
	\begin{equation}\label{eq:3.14}
	P(t,x,D_x)=\operatorname{diag}(\widetilde a_1,\ldots,\widetilde a_r)(t,x,D_x)+B(t,x,D_x),
	\end{equation}
	where $\widetilde a_i:=\operatorname{diag}(a_{k_i},\ldots,a_{k_{i+1}-1})$.
	
	Next, absorb the diagonal blocks of $B$ into the corresponding blocks of the principal part. Setting
	\[
	\widetilde A_i:=\widetilde a_i+\big(B_{\mu\nu}\big)_{k_i\le\mu,\nu\le k_{i+1}-1},
	\]
	we may assume (after this harmless redefinition) that the remaining matrix $B$ is block-off-diagonal, i.e.
	\begin{equation}\label{eq:3.15}
	B_{\mu\nu}=0\qquad\text{whenever }k_i\le\mu,\nu\le k_{i+1}-1.
	\end{equation}

	\smallskip
	\noindent\textit{Step 2: Factorization of the propagator.}
	Let $E(t,\tau)$ denote the block-diagonal propagator associated with the unperturbed system
	\[
	i\partial_t v-\widetilde A(t,x,D_x)v=0.
	\]
Thus $E(t,\tau)$ is a nondegenerate Fourier integral operator of order zero, with $E(t,\tau)$ and $E(t,\tau)^{-1}=E(\tau,t)$ bounded on $L^2$. No self-adjointness of $\widetilde A$ is required.
We factor the full propagator $U(t,\tau)$ for \eqref{eq:CP} as
	\[
	U(t,\tau)=E(t,\tau) V(t,\tau).
	\]
	Substituting into the equation for $U(t,\tau)$ and using that $E(t,\tau)$ solves the unperturbed problem, we obtain an evolution equation for $V(t,\tau)$:
	\begin{equation}\label{eq:3.16}
	\partial_t V(t,\tau) = Z(t,\tau) V(t,\tau),\qquad V(\tau,\tau)=I,
	\end{equation}
	where
	\[
	Z(t,\tau):=-i\,E(t,\tau)^{-1}\,B(t)\,E(t,\tau).
	\]
	For notational simplicity, we fix $\tau=0$ and write $V(t):=V(t,0)$, $Z(t):=Z(t,0)$.
	
	\smallskip
	\noindent\textit{Step 3: Application of Theorem~\ref{thm:3.1}.}
	We apply the abstract regularity principle to the Cauchy problem $V'(t)-Z(t)V(t)=0$ with initial datum $V(0)=I$.
	We choose the spaces
	\[
	W:=L^p\cap L^2_{\mathrm{comp}},\qquad
	W_1:=L^p_{\alpha}\cap L^2_{\mathrm{comp}},\qquad
	W_0:=H^s_{\mathrm{comp}},
	\]
	with $s>\frac{pn-2n}{2p}$ for $p>2$ and $s=0$ for $1<p\le 2$.	
	
	The hypotheses of Theorem~\ref{thm:3.1} are satisfied as follows.
	\begin{itemize}
		\item[(i)] \emph{Boundedness.} In general, $Z(t)$ is not a pseudo-differential operator. If $E_\mu$ and $E_\nu$ denote two diagonal block propagators, then an off-diagonal entry of $Z(t)$ has the form
		\[
		Z_{\mu\nu}(t)=-iE_\mu(0,t)B_{\mu\nu}(t)E_\nu(t,0),
		\qquad \mu\neq\nu.
		\]
		This is an order-zero Fourier integral operator associated with the composition of the canonical transformations generated by the two characteristic roots. The $L^2$ boundedness theorem for order-zero Fourier integral operators therefore shows that $Z(t)$ is uniformly bounded on $L^2$. Moreover, the Seeger--Sogge--Stein theorem \cite{SSS91} gives the local mapping property
		\[
		Z(t):L^p_\alpha\longrightarrow L^p_{\mathrm{loc}},
		\qquad \alpha=(n-1)\left|\frac1p-\frac12\right|.
		\]
		With the proper-support localisation used throughout, compactly supported inputs remain supported in a fixed compact set. Hence $Z(t)$ maps $W_1$ continuously into $W$, uniformly for $t\in[0,T]$.

		\item[(ii)] \emph{Stability under composition.} A finite composition
		$Z(t_1)\circ\cdots\circ Z(t_m)$ is again an order-zero Fourier integral operator, associated with the composed canonical transformation. Since the canonical relations are graphs, their composition remains a canonical graph. The Seeger--Sogge--Stein estimate therefore applies to the entire composition and yields a continuous map $W_1\to W$, uniformly for $(t_1,\ldots,t_m)\in[0,T]^m$. In particular, the loss of $\alpha$ derivatives is applied once to the composed Fourier integral operator, rather than once at each factor.

		\item[(iii)] \emph{Smoothing by time integration.} The Picard iterates $Q_\ell(t)$ associated with $Z$ are precisely those defined in \eqref{eq:Picard} after the substitution $V=E^{-1}U$. By Theorem~\ref{thm:smoothing},
		\[
		Q_\ell(t):L^2_{\mathrm{comp}}\longrightarrow
		H^{N(\ell)}_{\mathrm{loc}},
		\qquad N(\ell)\longrightarrow\infty.
		\]
		Choose $\ell$ so large that $N(\ell)\geq s$. Order-zero Fourier integral operators associated with canonical graphs are bounded on $H^s$, so $Z(t)Q_\ell(t)$ maps $W_1$ continuously into $H^s_{\mathrm{loc}}$. After the fixed proper-support localisation, this gives the required map into $W_0=H^s_{\mathrm{comp}}$.

		\item[(iv)] \emph{Regularity for zero Cauchy data.} Let $r(t)\in W_0$. Duhamel's formula gives
		\[
		v(t)=\int_0^t\widetilde V(t,\sigma)r(\sigma)\,d\sigma,
		\]
		where $\widetilde V(t,\sigma)$ is the propagator of $v'-Z(t)v=0$. Since $Z(t)$ is an order-zero Fourier integral operator and is uniformly bounded on $H^s$, the Picard series defining $\widetilde V(t,\sigma)$ converges in $\mathcal L(H^s)$. Thus $v(t)\in H^s_{\mathrm{loc}}$. The Sobolev embedding
		$H^s_{\mathrm{loc}}\hookrightarrow L^p_{\mathrm{loc}}$, valid for $s>n(\tfrac12-\tfrac1p)_+$, then gives $v(t)\in W$.
	\end{itemize}
	Theorem~\ref{thm:3.1} therefore yields $V(t):W_1\to W$ continuously, uniformly in $t\in[0,T]$.

	\smallskip
	\noindent\textit{Step 4: Conclusion for $u(t)$.}
	We do not claim that $V(t,0)$ is pseudo-differential or bounded on $L^p_\alpha$. Instead, set $Q_0=I$, choose $L$ so large that $Q_L$ has the Sobolev smoothing required in Step~3, and write
	\[
	u(t)=\sum_{k=0}^{L}E(t,0)Q_k(t)u_0
	      +E(t,0)R_{L+1}(t)u_0,
	\qquad
	R_{L+1}(t):=\sum_{k=L+1}^{\infty}Q_k(t).
	\]
	For each fixed $k$, the operator $E(t,0)Q_k(t)$ is a time integral of a finite sum of order-zero Fourier integral operators associated with canonical graphs. The Seeger--Sogge--Stein theorem therefore gives
	\[
	E(t,0)Q_k(t):L^p_\alpha\longrightarrow L^p_{\mathrm{loc}}.
	\]
	The remainder satisfies
	\[
	R_{L+1}'(t)=Z(t)R_{L+1}(t)+Z(t)Q_L(t),
	\qquad R_{L+1}(0)=0.
	\]
	By Step~3 and Duhamel's formula, $R_{L+1}(t)u_0\in H^s_{\mathrm{loc}}$; since $E(t,0)$ is an order-zero Fourier integral operator, it is bounded on local Sobolev spaces. Sobolev embedding therefore controls the remainder in $L^p_{\mathrm{loc}}$. Combining the finite Picard sum with this remainder estimate yields
	\[
	\|u(t,\cdot)\|_{L^p_{\mathrm{loc}}}
	\leq C_T\|u_0\|_{L^p_\alpha},
	\]
	uniformly for $t\in[0,T]$. This is exactly~\eqref{eq:2.4}.
\qed

\subsection{Proof of Theorem~\ref{thm:weak-eq}}

We reduce the scalar equation to a first-order system and apply Theorem~\ref{thm:Lp-reg}.

\smallskip
\noindent\textit{Step 1: A first-order system for derivatives of $u$.}
Recall that $\partial_j:=D_t+\lambda_j(t,x,D_x)$ for $j=1,\ldots,m$, so that
\[
L=\partial_1\cdots\partial_m+\sum_{k<m} b_{j_1,\ldots,j_k}\,\partial_{j_1}\cdots\partial_{j_k}+c.
\]
Following Morimoto \cite{Mor79}, let $U$ be the vector of all iterated derivatives of length at most $m-1$,
\[
U:=\big(\partial^J u\big)_{|J|\le m-1},
\]
where $J=(j_1,\ldots,j_k)$ with $1\le j_\ell\le m$ and $\partial^J:=\partial_{j_1}\cdots\partial_{j_k}$ (the components are ordered arbitrarily but fixed). Then $U$ satisfies a first-order system
\begin{equation}\label{eq:3.18}
D_t U + A(t,x,D_x)U + B(t,x,D_x)U = F,
\qquad U(0,\cdot)=G,
\end{equation}
with $A$ diagonal with entries $\lambda_j$, $B$ of order $0$, and $F$ depending linearly on $Lu$. Under the hypothesis $Lu=0$ of Theorem~\ref{thm:weak-eq}, we have $F\equiv 0$.

\smallskip
\noindent\textit{Step 2: Identification of the initial vector $G$.}
For each multi-index $J$ with $|J|=k\le m-1$, expand
\begin{equation}\label{eq:thm23_expansion}
\partial^J = D_t^k + \sum_{j=0}^{k-1} c_j^J(t,x,D_x)\,D_t^j,\qquad c_j^J\in\Psi^{k-j}.
\end{equation}
Evaluating at $t=0$ and using $D_t^j u(0,\cdot)=g_j$ gives
\[
(\partial^J u)(0,\cdot)=g_k+\sum_{j=0}^{k-1} c_j^J(0,x,D_x)\,g_j,
\]
so each component of $G$ is a finite sum of pseudo-differential operators of order $k-j$ applied to $g_j$. Since the $g_j$ are compactly supported, the mapping properties of $\Psi^{k-j}$ yield
\[
\|G\|_{L^p_{\alpha}} \le C\sum_{j=0}^{m-1}\|g_j\|_{L^p_{\alpha+m-1-j}}.
\]

\smallskip
\noindent\textit{Step 3: Application of Theorem~\ref{thm:Lp-reg}.}
Under Condition~C for the principal symbols of the $\lambda_j$, the system \eqref{eq:3.18} has diagonal principal part and an order-zero perturbation, fitting the framework of Theorem~\ref{thm:Lp-reg}. Applying the theorem to \eqref{eq:3.18} with $F\equiv 0$ gives
\[
\sup_{t\in[0,T]}\|U(t,\cdot)\|_{L^p_{\mathrm{loc}}}
\le C\,\|G\|_{L^p_{\alpha}}
\le C\sum_{j=0}^{m-1}\|g_j\|_{L^p_{\alpha+m-1-j}}.
\]
Since $U$ contains every $\partial^J u$ with $|J|\le m-1$, this yields the stated estimate for each such derivative. Applying $(I-\Delta_x)^{s/2}$ (equivalently, working in $L^p_s$ throughout) extends the bound to $(L^p_s)_{\mathrm{loc}}$, completing the proof. \qed

\section{Representation formula for non-diagonalisable systems} \label{sec:triangular-representation}
In this section, we remove the assumption in \cite{GJR20} that the principal part is time-independent, and extend the representation formula obtained there for non-diagonalisable systems with variable multiplicities to time-dependent principal parts. This completes the study of hyperbolic systems with non-diagonalisable principal part and variable multiplicities initiated in \cite{GJR18, GJR20} and continued in \cite{GS24}.

Let $n\ge 1$, $m\ge 2$, $T>0$, and write $D_t=-i\partial_t$, $D_x=-i\partial_x$.
We consider the Cauchy problem
\begin{equation}\label{rep:triangular-system}
\begin{cases}
D_t u = A(t,x,D_x)\,u + B(t,x,D_x)\,u + f(t,x), & (t,x)\in[0,T]\times\mathbb R^n,\\
u(0,x)=u_0(x), & x\in\mathbb R^n,
\end{cases}
\end{equation}
where $u=(u_1,\dots,u_m)^{T}$ is the unknown, $f$ and $u_0$ are given
$\mathbb C^m$-valued data, and $A=(a_{ij})$, $B=(b_{ij})$ are $m\times m$
matrices of pseudo-differential operators with
$a_{ij}\in C^\infty([0,T];\Psi^1(\mathbb R^n))$ and
$b_{ij}\in C^\infty([0,T];\Psi^0(\mathbb R^n))$.
We denote by $a_j(t,x,\xi)$ the principal symbol of $a_{jj}$.

\begin{assumption}\label{rep:assumptions}
\leavevmode
\begin{enumerate}
  \item[(H1)]\label{H:triangular} \emph{(Triangular structure)} $A$ is upper triangular,
  so $A=\Lambda+N$ with $\Lambda=\operatorname{diag}(a_{11},\dots,a_{mm})$
  and $N$ strictly upper triangular.
  \item[(H2)]\label{H:hyperbolic} \emph{(Hyperbolicity)} Each $a_j$ is real-valued and
  positively homogeneous of degree one in $\xi$.
  \item[(H3)]\label{H:levi} \emph{(Levi-type condition)}
  $b_{ij}\in C^\infty([0,T];\Psi^{j-i}(\mathbb R^n))$ for $i>j$.
  \item[(H4)]\label{H:C} The principal symbols $a_1,\dots,a_m$ satisfy Condition C.
\end{enumerate}
\end{assumption}

Under these assumptions, well-posedness of \eqref{rep:triangular-system} in anisotropic Sobolev spaces was proved in \cite{GJR18}, and a representation formula was obtained in \cite{GJR20} when the principal part is time-independent.

The triangular structure (H1) suggests solving the system
through its diagonal entries. For $j=1,\dots,m$, let $E_j(t,\tau)$ denote the
solution operator of the scalar Cauchy problem
\[
\begin{cases}
D_t E_j(t,\tau)=\bigl(a_{jj}+b_{jj}\bigr)(t,x,D_x)\,E_j(t,\tau),\\
E_j(\tau,\tau)=I,
\end{cases}
\]
which is well posed by (H2). By uniqueness,
$E_j(t,s)E_j(s,\tau)=E_j(t,\tau)$ for all $t,s,\tau\in[0,T]$; in particular,
$E_j(t,\tau)^{-1}=E_j(\tau,t)$. We define
\begin{equation}\label{rep:scalar-operators}
(G_j^0 g)(t)=E_j(t,0)\,g,
\qquad
(G_j f)(t)=i\int_0^t E_j(t,\tau)\,f(\tau)\,d\tau .
\end{equation}
Then $w=G_j^0 g+G_j f$ is the solution of
$D_t w=(a_{jj}+b_{jj})w+f$, $w(0)=g$.

\begin{theorem}\label{rep:thm-main}
Let Assumption~\ref{rep:assumptions} hold and let $s\in\mathbb R$. Suppose that
$u_0=(u_1^0,\dots,u_m^0)^T$ and $f=(f_1,\dots,f_m)^T$ satisfy
\[
u_l^0\in H^{s+l-1}(\mathbb R^n),\qquad
f_l\in C([0,T];H^{s+l-1}(\mathbb R^n)),\qquad l=1,\dots,m.
\]
Then the following hold.
\begin{enumerate}
  \item Problem \eqref{rep:triangular-system} has a unique solution $u$ with
  $u_j\in C([0,T];H^{s+j-1}(\mathbb R^n))$ for $j=1,\dots,m$.
  \item For every $N\in\mathbb N$, the components of $u$ are given by
  \begin{equation}\label{rep:representation}
  u_j(t)=\sum_{l=1}^{m}\Bigl[\bigl(\mathcal H_{j,l}(t)+R_{j,l}(t)\bigr)u_l^0
  +\bigl(\mathcal K_{j,l}(t)+S_{j,l}(t)\bigr)f_l\Bigr],
  \qquad j=1,\dots,m,
  \end{equation}
  where $\mathcal H_{j,l}$ and $\mathcal K_{j,l}$ are integrated Fourier
  integral operators of order $l-j$, with
  \[
  \mathcal H_{j,l}\in\mathcal L\bigl(H^{s+l-1},\,C([0,T];H^{s+j-1})\bigr),\qquad
  \mathcal K_{j,l}\in\mathcal L\bigl(C([0,T];H^{s+l-1}),\,C([0,T];H^{s+j-1})\bigr),
  \]
  and the remainders satisfy
  \[
  R_{j,l}\in\mathcal L\bigl(H^{s+l-1},\,C([0,T];H^{s+j-1+N})\bigr),\qquad
  S_{j,l}\in\mathcal L\bigl(C([0,T];H^{s+l-1}),\,C([0,T];H^{s+j-1+N})\bigr).
  \]
\end{enumerate}
\end{theorem}
The proof follows the lines of \cite[Theorem~1.3]{GJR20}. The only point where time-independence of the principal part was used there is that the scalar propagators satisfy $E_j(t,\tau)=E_j(t-\tau)$. Here $E_j(t,\tau)$ depends on $t$ and $\tau$ separately, and we show that the argument only requires the group property $E_j(t,s)E_j(s,\tau)=E_j(t,\tau)$. We give the details for $m=2$, which already contains all new ingredients, and then indicate the changes for general $m$.

\subsection{The case $m=2$}
As a model problem, we consider \eqref{rep:triangular-system} with $m=2$.
It contains all new ingredients, and the case $m\ge 3$ follows by the same argument as in \cite{GJR20}.
\subsubsection*{Reduction to scalar equations.}
For $m=2$ the system reads
\[
D_tu_1=(a_{11}+b_{11})u_1+(a_{12}+b_{12})u_2+f_1,\qquad
D_tu_2=(a_{22}+b_{22})u_2+b_{21}u_1+f_2,
\]
where $a_{12}$ has order one and $b_{21}\in\Psi^{-1}$ by (H3). Treating the coupling terms as sources, \eqref{rep:scalar-operators} gives
\begin{align}
u_1&=G_1^0u_1^0+G_1f_1+G_1(a_{12}+b_{12})u_2,\label{rep:m2-u1}\\
u_2&=G_2^0u_2^0+G_2f_2+G_2b_{21}u_1.\label{rep:m2-u2}
\end{align}
Substituting each equation into the other, we obtain
\begin{align*}
L_1u_1&=G_1^0u_1^0+G_1f_1+G_1(a_{12}+b_{12})\bigl(G_2^0u_2^0+G_2f_2\bigr),\\
L_2u_2&=G_2^0u_2^0+G_2f_2+G_2b_{21}\bigl(G_1^0u_1^0+G_1f_1\bigr),
\end{align*}
where $L_j=I-\mathcal G_j$ and
\[
\mathcal G_1=G_1(a_{12}+b_{12})G_2b_{21},\qquad
\mathcal G_2=G_2b_{21}G_1(a_{12}+b_{12}).
\]

\subsubsection*{Inverting $L_j$.}
Since $b_{21}$ has order $-1$, each $\mathcal G_j$ has order zero. Since
$G_j$ contains a time integral, $\mathcal G_j$ is of Volterra type:
\[
\|(\mathcal G_j^kw)(t)\|_{H^s}\le\frac{(Ct)^k}{k!}\,
\sup_{0\le\tau\le t}\|w(\tau)\|_{H^s},\qquad k\ge0 .
\]
Hence $L_j^{-1}=\sum_{k\ge0}\mathcal G_j^k$ converges in
$\mathcal L\bigl(C([0,T];H^s)\bigr)$ for every $s\in\mathbb R$, and
\begin{equation}\label{rep:m2-solution}
\begin{aligned}
u_1&=L_1^{-1}\Bigl[G_1^0u_1^0+G_1f_1+G_1(a_{12}+b_{12})\bigl(G_2^0u_2^0+G_2f_2\bigr)\Bigr],\\
u_2&=L_2^{-1}\Bigl[G_2^0u_2^0+G_2f_2+G_2b_{21}\bigl(G_1^0u_1^0+G_1f_1\bigr)\Bigr].
\end{aligned}
\end{equation}

\begin{remark}\label{rep:rem-picard}
Formula \eqref{rep:m2-solution} agrees with the Picard series approach of
\cite{GJR20}. We show this for the homogeneous part, that is, for the
fundamental solution $U(t,\tau)$ of
\[
D_tU(t,\tau)=\bigl(\widetilde\Lambda(t)+C(t)\bigr)U(t,\tau),\qquad U(\tau,\tau)=I,
\]
where
\[
\widetilde\Lambda=\begin{pmatrix}a_{11}+b_{11}&0\\0&a_{22}+b_{22}\end{pmatrix},
\qquad
C=\begin{pmatrix}0&a_{12}+b_{12}\\ b_{21}&0\end{pmatrix}.
\]

\smallskip
\noindent\emph{Step 1 (Factorisation).}
Let $E(t,\tau)=\operatorname{diag}\bigl(E_1(t,\tau),E_2(t,\tau)\bigr)$ and write
$U=EV$. Since $D_tE=\widetilde\Lambda E$, we get $E\,D_tV=C\,E\,V$, and since
$E(t,\tau)^{-1}=E(\tau,t)$,
\[
D_tV(t,\tau)=Z(t,\tau)V(t,\tau),\qquad V(\tau,\tau)=I,\qquad
Z(t,\tau)=E(\tau,t)\,C(t)\,E(t,\tau).
\]

\smallskip
\noindent\emph{Step 2 (Picard series).}
Writing the equation for $V$ in integral form gives
$V=\sum_{k\ge0}V_k$ with
\[
V_0=I,\qquad
V_k(t,\tau)=i\int_\tau^tZ(s,\tau)\,V_{k-1}(s,\tau)\,ds,\quad k\ge1 .
\]
Since $Z$ is off-diagonal, $V_k$ is diagonal for even $k$ and off-diagonal
for odd $k$.

\smallskip
\noindent\emph{Step 3 (Recursion in terms of $G_j$).}
From now on $\tau=0$. Set $H_k(t)=E(t,0)V_k(t,0)$, so $H_0=E(\cdot,0)$.
The group law gives
\[
E(t,0)Z(s,0)=E(t,0)E(0,s)\,C(s)\,E(s,0)=E(t,s)\,C(s)\,E(s,0),
\]
and therefore, for $k\ge1$,
\[
H_k(t)=i\int_0^tE(t,s)\,C(s)\,H_{k-1}(s)\,ds=(\mathcal TH_{k-1})(t),
\qquad
\mathcal T=\begin{pmatrix}0&G_1(a_{12}+b_{12})\\ G_2b_{21}&0\end{pmatrix},
\]
where $\mathcal T$ acts on each column. Hence $H_k=\mathcal T^kE(\cdot,0)$.

\smallskip
\noindent\emph{Step 4 (Powers of $\mathcal T$).}
Direct multiplication gives $\mathcal T^2=\operatorname{diag}(\mathcal G_1,\mathcal G_2)$,
so for every $k\ge0$,
\[
\mathcal T^{2k}=\begin{pmatrix}\mathcal G_1^k&0\\0&\mathcal G_2^k\end{pmatrix},
\qquad
\mathcal T^{2k+1}=\begin{pmatrix}0&\mathcal G_1^kG_1(a_{12}+b_{12})\\
\mathcal G_2^kG_2b_{21}&0\end{pmatrix}.
\]

\smallskip
\noindent\emph{Step 5 (Summation).}
Since $U(t,0)=\sum_{k\ge0}H_k(t)$, we split the sum into even and odd $k$.
The Neumann series $\sum_{k\ge0}\mathcal G_j^k=L_j^{-1}$ converges, which
justifies the rearrangement, and we obtain
\[
U(t,0)=
\begin{pmatrix}
L_1^{-1}E_1(\cdot,0) & L_1^{-1}G_1(a_{12}+b_{12})E_2(\cdot,0)\\
L_2^{-1}G_2b_{21}E_1(\cdot,0) & L_2^{-1}E_2(\cdot,0)
\end{pmatrix}(t).
\]
Applied to $u_0$, and recalling $G_j^0=E_j(\cdot,0)$, this is exactly the
$u_0$ part of \eqref{rep:m2-solution}.

\smallskip
In \cite{GJR20} the principal part is time-independent and
$E_j(t,s)=E_j(t-s)$. The argument above uses only the group law
$E_j(t,s)E_j(s,\tau)=E_j(t,\tau)$, so it applies without change to
time-dependent principal parts.
\end{remark}

\subsubsection*{Representation formula.}
Fix $N\in\mathbb N$. By Theorem~\ref{thm:smoothing} and Condition~C, there is
$M=M(N)$ such that
\[
L_1^{-1}=\sum_{k=0}^{M}\mathcal G_1^k+R_1,\qquad
R_1\in\mathcal L\bigl(C([0,T];H^s),\,C([0,T];H^{s+N})\bigr)
\ \text{for all } s\in\mathbb R .
\]
Inserting this into \eqref{rep:m2-solution} gives \eqref{rep:representation}
for $j=1$, with
\begin{align*}
\mathcal H_{1,1}&=\sum_{k=0}^{M}\mathcal G_1^kG_1^0, &
\mathcal H_{1,2}&=\sum_{k=0}^{M}\mathcal G_1^kG_1(a_{12}+b_{12})G_2^0,\\
\mathcal K_{1,1}&=\sum_{k=0}^{M}\mathcal G_1^kG_1, &
\mathcal K_{1,2}&=\sum_{k=0}^{M}\mathcal G_1^kG_1(a_{12}+b_{12})G_2,
\end{align*}
and remainders $R_{1,1}=R_1G_1^0$, $R_{1,2}=R_1G_1(a_{12}+b_{12})G_2^0$,
$S_{1,1}=R_1G_1$, $S_{1,2}=R_1G_1(a_{12}+b_{12})G_2$. Since $a_{12}+b_{12}$ has order one and all other factors have order zero, $\mathcal H_{1,l}$ and $\mathcal K_{1,l}$ have order $l-1$, and
$R_{1,l}$, $S_{1,l}$ map $H^{s+l-1}$ and $C([0,T];H^{s+l-1})$, respectively, into $C([0,T];H^{s+N})$, as stated in Theorem~\ref{rep:thm-main}. By the composition results of \cite{GJR20}, $\mathcal H_{1,l}$ and
$\mathcal K_{1,l}$ are integrated Fourier integral operators. The same argument applied to $u_2$, with $\mathcal G_2$ in place of $\mathcal G_1$, gives \eqref{rep:representation} for $j=2$; for instance,
$\mathcal H_{2,1}=\sum_{k=0}^{M}\mathcal G_2^kG_2b_{21}G_1^0$ has order $-1$.

\appendix
\section{Auxiliary Lemmas}\label{sec:auxiliary}
We collect here five technical results used in the proofs of Theorems~\ref{thm:smoothing} and~\ref{thm:Lp-reg}: a lifting lemma ensuring global diagonalizability of the principal symbol (\S\ref{app:diag}); Egorov's theorem on conjugation by Fourier integral propagators (\S\ref{app:egorov}); an $L^2\to H^s$ estimate for FIOs (\S\ref{app:FIO}); a measure estimate for sublevel sets of smooth functions with finite-order zeros (\S\ref{app:meas}); and an interpolation lemma (\S\ref{app:interp}). Lemmas~\ref{lem:diag}, \ref{thm:egorov}, and~\ref{lem:interpolation} are recalled from the literature; Lemma~\ref{lem:FIO-est} adapts \cite[Lemma~5.2]{KR2007} to $\mathbb{R}^n$; Proposition~\ref{prop:measure} is the standard measure estimate, included for completeness.

\subsection{Lifting and global diagonalization}\label{app:diag}

To make global use of the diagonalization of the principal symbol, one must first ensure that its characteristic roots admit smooth global branches.n general, these roots are defined only locally, since analytic continuation around closed loops in $S^*X$ may induce a permutation among them. The following lemma, due to \cite{KR2007}, resolves this by passing to a finite cover $\widetilde X \to X$ on which the eigenvalue branches become single-valued. The dimensions of $X$ and $\widetilde X$ coincide, so estimates and mapping properties on $\widetilde X$ descend to $X$ at the end of the argument.  The construction is analogous to that given in the self-adjoint setting by \cite{Roz80}.

The restriction $n\ge 3$ enters here: in lower dimensions the homotopy-theoretic step below fails, and one must instead either assume that the principal symbol is locally smoothly diagonalizable on $X$ itself or work microlocally throughout. The main theorems of this paper (Theorems~\ref{thm:smoothing}--\ref{thm:weak-eq}) are stated on $\mathbb{R}^n$ and use the local lifting argument of Remark~\ref{rem:covering}, which holds in every dimension $n\ge 1$; consequently they are not subject to the $n\ge 3$ restriction. Lemma~\ref{lem:diag} is recorded here only for the global manifold setting (where it is the natural starting point for spectral asymptotics).

\begin{lemma}[\cite{KR2007}]\label{lem:diag}
	Let $X$ be a smooth compact manifold without boundary, $\dim X=n\ge 3$,
	and let $P(x,D)$ be a first-order pseudo-differential operator acting
	on sections of an $m$-dimensional vector bundle $E\to X$. Denote by
	$L^2(E)$ the space of half-density sections, and by $E'$ the lift of
	$E$ to $T^*X$. Assume that the principal symbol $A(x,\xi)$ is
	microlocally smoothly diagonalizable: in any cone
	$\Lambda\subset T^*X$ on which $E'|_\Lambda\simeq\Lambda\times\mathbb C^m$,
	the matrix $A(x,\xi)|_\Lambda$ admits $m$ smooth real eigenvalues
	$a_j(x,\xi)$ with one-dimensional eigenspaces $V_j(x,\xi)$, and these
	diagonalizations agree on overlaps.
	
	Then there exists a finite cover $p:\widetilde X\to X$ such that, on
	the lift $\widetilde E$ of $E$ to $T^*\widetilde X$, the eigenvalues
	$a_j$ and eigenspaces $V_j$ admit smooth global branches. The space
	$L^2(\widetilde E)$ decomposes as a direct sum of $m$ subspaces with
	respect to which $P(\widetilde x,\widetilde D)$ is represented by a
	matrix of pseudo-differential operators whose principal symbol is the
	diagonal matrix $\operatorname{diag}\big(a_j(\widetilde x,\widetilde\xi)\big)$.
\end{lemma}

\begin{remark}\label{rem:5.1-time-dependent}
	The same conclusion holds for a smooth family $P(t,x,D_x)\in C^{\infty}([0,T],\Psi^1)$ whose principal symbol
	$A(t,x,\xi)$ is microlocally smoothly diagonalizable jointly in $(t,x,\xi)\in[0,T]\times T^*X\setminus 0$. Indeed, for each fixed $t$ the analytic continuation of the eigenvalues along closed paths in $S^*X$ defines a homomorphism
	\[
	\rho_t:\pi_1(X)\longrightarrow S_m,
	\]
	and since $\rho_t$ takes values in a discrete group and depends
	continuously on $t$, it is constant on the connected interval
	$[0,T]$. Consequently, the cover $p:\widetilde X\to X$ produced by
	Lemma~\ref{lem:diag} is \emph{independent of $t$}. The orthogonal
	projectors $p_j(t,\widetilde x,\widetilde D)$ depend smoothly on $t$,
	and the resulting diagonalization
	\[
	P(t,\widetilde x,\widetilde D)
	=\sum_{j,k} p_j(t)\,P(t)\,p_k(t)
	\]
	yields a smooth family in $C^{\infty}([0,T],\Psi^1)$ with principal symbol
	$\operatorname{diag}\big(a_1(t,\widetilde x,\widetilde\xi),\ldots,a_m(t,\widetilde x,\widetilde\xi)\big)$.
	
	When eigenvalues coincide along $[0,T]$ but separate into smooth
	branches of constant multiplicity, the same argument yields a
	block-diagonal representation with smooth principal blocks. This is
	the form in which the lemma enters the proof of
	Theorem~\ref{thm:Lp-reg}.
\end{remark}

\subsection{Conjugation by propagators}\label{app:egorov}

We will make essential use of the following lemma, which describes how pseudo-differential operators transform under conjugation by propagators.

\begin{lemma}[Egorov's theorem \cite{T81}]\label{thm:egorov}
Let $\widehat A_j(t,x,D_x) \in C^\infty([0,T],\Psi^1)$ have real principal
symbol $a_j(t,x,\xi)$, and let $E_j(t,s)$ be its propagator, i.e.\ 
the solution of
\[
i\partial_t E_j(t,s) = \widehat A_j(t,x,D_x)\, E_j(t,s), \qquad E_j(s,s) = I,
\qquad t,s\in[0,T].
\]
No self-adjointness is assumed: the conjugation uses
$E_j(s,t)=E_j(t,s)^{-1}$, not the adjoint.
Let $R(t)\in C^\infty([0,T],\Psi^\mu)$ have principal symbol 
$r(t,x,\xi)$. Then the conjugated family
\[
\widetilde{R}(t,s) := E_j(t,s)\, R(s)\, E_j(s,t)
\]
belongs to $\Psi^\mu$ for every $t,s\in[0,T]$, and its principal 
symbol is given by
\begin{equation}\label{eq:egorov-symbol}
\sigma_\mu\big(\widetilde{R}(t,s)\big)(x,\xi) 
= r\bigl(s,\, \Phi_j^{s,t}(x,\xi)\bigr) \pmod{S^{\mu-1}},
\end{equation}
where $\Phi_j^{s,t}\colon T^*\mathbb{R}^n \to T^*\mathbb{R}^n$ 
denotes the Hamiltonian flow generated by $a_j$ from time $t$ to 
time $s$.
\end{lemma}

\subsection{Fourier integral operator estimates}\label{app:FIO}
The following lemma gives $L^2 \to H^s$ estimates for Fourier integral  operators in terms of $C^k$ norms of their amplitudes. It is a variant of \cite[Lemma~5.2]{KR2007}, adapted to the setting of $\mathbb{R}^n$.

\begin{lemma}\label{lem:FIO-est}
	Let $X \subset \mathbb{R}^n$ be open, let $q \in \mathbb{Z}$, and 
	let $T$ be a Fourier integral operator of the form
	\begin{equation}\label{eq:FIO-form}
	(Tu)(x) = \int_{\mathbb{R}^n} \int_{\mathbb{R}^n} 
	e^{i\varphi(x,y,\xi)}\, a(x,y,\xi)\, u(y)\, dy\, d\xi,
	\end{equation}
	where the amplitude $a \in S^q_{1,0}$ has compact support in 
	$(x,y)\in X\times X$ uniformly in $\xi$, and the phase satisfies the 
	non-degeneracy condition $\det\partial^2_{x,\xi}\varphi\neq 0$ for $\xi\ne 0$. 
	Then $T$ extends to a bounded operator
	\(
	T\colon L^2(X) \to H^{-q-n-1}(X),
	\)
	with
	\begin{equation}\label{eq:FIO-L2Hs}
	\|T\|_{L^2 \to H^{-q-n-1}} 
	\leq C\, \big\|\, a\,\langle\xi\rangle^{-q}\,\big\|_{C^{|q+n+1|}},
	\end{equation}
	where $\langle\xi\rangle := (1+|\xi|^2)^{1/2}$, and the constant 
	$C$ depends only on $n$ and $q$.
	
	Suppose in addition that the amplitude is supported in a $\xi$-conical 
	subset of $X\times X\times(\mathbb{R}^n\setminus 0)$,
	\begin{equation}\label{eq:conical-set}
	\operatorname{supp} a \;\subset\; \Xi_1 
	:= \left\{(x,y,\xi)\colon\; x \in X,\; 
	\Big(y,\tfrac{\xi}{|\xi|}\Big) \in \Xi,\; \xi \neq 0\right\},
	\end{equation}
	where $\Xi\subset S^*X$ is a measurable set, viewed in coordinates 
	$(y,\eta)$ with $|\eta|=1$. Then, for every $\delta > 0$,
	\begin{equation}\label{eq:FIO-meas}
	\|T\|_{L^2 \to H^{-q-3n/2-1-\delta}} 
	\leq C\, \mathrm{meas}(\Xi)\, 
	\big\|\, a\,\langle\xi\rangle^{-q}\,\big\|_{C^{|q+n+1|}},
	\end{equation}
	where $\mathrm{meas}$ denotes the canonical induced measure on 
	$S^*X$ and $C$ depends on $n$, $q$, $\delta$, and on the size of 
	the support of $a$ in $(x,y)$.
\end{lemma}

\begin{remark}
	The additional loss of $n/2+\delta$ derivatives in 
	\eqref{eq:FIO-meas} compared with \eqref{eq:FIO-L2Hs} is the price 
	for replacing the $L^\infty$-control of the amplitude by the 
	measure of its $\xi$-support: it arises from a Cauchy--Schwarz 
	estimate over the slice $\Xi\subset S^*X$.
\end{remark}

\subsection{A measure estimate for finite-order zeros}\label{app:meas}
\label{subsec:measure}

The following proposition controls the measure of the sublevel set on which a smooth function with zeros of uniformly bounded order is small. It enters the proof of Theorem~\ref{thm:smoothing} (Stage~III) where it is used to estimate $\mathrm{meas}(\Xi)$ under Condition~C.

\begin{proposition}\label{prop:measure}
	Let $W\subset\mathbb{R}^n$ be compact, $T\in(0,\infty)$, and $M\in\mathbb{N}$. Let $f\colon[0,T]\times W\to\mathbb{R}$ be continuous in $w$ and smooth in $t$, and assume that for every $(t,w)\in[0,T]\times W$ there exists $k\in\{1,\ldots,M\}$ with
	\[
	|\partial_t^k f(t,w)|\ne 0,
	\]
	i.e., $f(\cdot,w)$ vanishes to order at most $M$ at each of its zeros, uniformly in $w\in W$. Then there exist constants $C>0$ and $\varepsilon_0>0$ such that, for every $0<\varepsilon<\varepsilon_0$,
	\begin{equation}\label{eq:measure-est}
	\sup_{w\in W}\,\mathrm{meas}\bigl\{t\in[0,T]:|f(t,w)|\leq\varepsilon\bigr\}\;\leq\;C\,\varepsilon^{1/2M}.
	\end{equation}
\end{proposition}

\subsection{Interpolation}\label{app:interp}

The final ingredient is the following interpolation lemma 
(\cite[Lemma~5.6]{KR2007}), which we use at the end of Stage~III to 
combine the estimates for $R_1$, $R_2$, and $R_3$ into the smoothing 
estimate for $Q_\ell$.

\begin{lemma}[Interpolation, \cite{KR2007}]\label{lem:interpolation}
	Let $p<r$, $\sigma,\tau>0$, and $u\in H^p(\mathbb{R}^n)$. Suppose 
	that for every sufficiently small $\varepsilon>0$ there exists a 
	decomposition
	\[
	u = u_1^{(\varepsilon)} + u_2^{(\varepsilon)}
	\]
	with
	\begin{equation}\label{eq:interp-hyp}
	\big\|u_1^{(\varepsilon)}\big\|_{H^r} \leq C\varepsilon^{-\tau},
	\qquad 
	\big\|u_2^{(\varepsilon)}\big\|_{H^p} \leq C\varepsilon^{\sigma},
	\end{equation}
	for a constant $C$ independent of $\varepsilon$. Then 
	$u\in H^q(\mathbb{R}^n)$ for every
	\begin{equation}\label{eq:interp-q}
	q < \frac{p\tau + r\sigma}{\tau + \sigma}.
	\end{equation}
\end{lemma}

\begin{remark}
	The exponent in \eqref{eq:interp-q} is the weighted average of $p$ 
	and $r$ with weights $\tau$ and $\sigma$; it interpolates between 
	the trivial regularity $u\in H^p$ (when $\sigma=0$) and the full 
	regularity $u\in H^r$ (when $\tau=0$).
\end{remark}

\subsection*{Acknowledgements}
M.~Ruzhansky and B.~Sabitbek were supported by the Methusalem programme of the Ghent University Special Research Fund (BOF) (Grant Nos.\ BOFMET2021000601 and 01M01021). M.~Ruzhansky was also supported by EPSRC Grant UKRI3645 and by FWO Senior Research Grant G022821N.

\subsection*{Use of AI tools}
During the preparation of this work, the authors used Claude (Anthropic) to improve the clarity and language of the manuscript and to assist in checking computations. All mathematical results, proofs and conclusions were developed by the authors, who reviewed and edited all content and take full responsibility for the publication.


\end{document}